\documentclass[12pt,reqno]{amsart}
\usepackage{amscd,amssymb,euscript,amstext , amsthm}
\title{Noncommutative Poincar\'e duality for boundary actions of hyperbolic groups}
\author{Heath Emerson}
\date{November 21, 2002}
\parskip=\baselineskip
\newcounter{axiom}[axiom]
\newcommand{\ebar}{\underline{E}\Gamma}

\theoremstyle{plain}
\newtheorem{thm}{Theorem}
\newtheorem{lemma}[thm]{Lemma}

\newtheorem{prop}[thm]{Proposition}
\newtheorem{cor}[thm]{Corollary}
\newtheorem{defn}[thm]{Definition}

\newtheorem{rmk}[thm]{Remark}

\newtheorem{note}[thm]{Note}
\newcommand{\C}{\mathbb C}

\newcommand{\Z}{\mathbb Z}
\newcommand{\R}{\mathbb R}

\newcommand{\F}{\mathbb F}

\newcommand{\Dudelta}{\hat{\Delta}}

\newcommand{\bgamma}{\partial \Gamma}

\newcommand{\elambda}{\underline{E}\Lambda}

\begin{document}
\maketitle

\begin{abstract}

For a large class of word hyperbolic groups $\Gamma$ the
cross product $C^{*}$-algebras $C(\bgamma) \rtimes \Gamma$, where
$\bgamma$ denotes the Gromov boundary of $\Gamma$ satisfy
Poincar\'e duality in $K$-theory. This class strictly contains fundamental groups of compact, negatively
curved manifolds. We discuss the general notion of Poincar\'e duality
for $C^{*}$-algebras, construct the fundamental classes for the
aforementioned algebras, and prove that $KK$-products with these classes induce inverse isomorphisms. The Baum-Connes Conjecture for
amenable groupoids is used in a crucial way.

\end{abstract}

\section{Introduction}

It is well known that if $M^{n}$ is a compact $n$-dimensional spin$^{c}$-manifold, the
$C^{*}$-algebra $C(M^{n})$ of continuous functions on $M^{n}$ exhibits
Poincar\'e duality in $K$-theory. Specifically, the class $[D] \in
K_{n}(M)$ of the Dirac
operator on $M^{n}$ induces by cap product an isomorphism $K^{*}(M^{n})
\to K_{* + n}(M^{n})$.  It is natural to ask whether
there are noncommutative $C^{*}$-algebras exhibiting the same
phenomenon. In \cite{Co2} A. Connes introduced the appropriate formalism for this
question, defining the analog for $C^{*}$-algebras of
Spanier-Whitehead duality for finite complexes. Two $C^{*}$-algebras
$A$ and $B$ shall be said to be dual if there exists a class $\Delta$
in the $K$-homology of $A \otimes B$, and a class $\Dudelta$ in the $K$-theory
of $A \otimes B$ such that $\Dudelta \otimes_{B} \Delta = 1_{A}$ and
$\Dudelta \otimes_{A} \Delta = 1_{B}$. If $A$ and $B$ are dual, cap
product with $\Delta $ induces an isomorphism $K_{*}(A) \to
K^{*}(B)$. A special case is
where $B = A^{\rm op}$, which we term Poincar\'e duality, while a 
$C^{*}$-algebra satisfying Poincar\'e duality we shall call in this paper a
Poincar\'e duality algebra. Known $commutative$ examples of Poincar\'e
duality algebras are given by continuous functions on spaces homotopy
equivalent to one of the aforementioned $M^{n}$ above; it is unknown
to the author whether there are other commutative examples. The first
nontrivial example of a noncommutative Poincar\'e duality algebra was
given by Connes
(see \cite{Co1}) in the form of the irrational rotation algebra
$A_{\theta}$.  In this
paper we shall prove that if $\Gamma$ is a hyperbolic group satisfying
a certain mild symmetry property, and $\bgamma$ is its Gromov boundary,
then the cross product $C(\bgamma) \rtimes \Gamma$ is a Poincar\'e
duality algebra.

Examples of pairs of algebras $A$ and $B$ dual in the above sense were given by Kaminker and Putnam (see \cite{KP}); the pairs
were $O_{M}$ and $O_{M^{t}}$ respectively, where for a square
$0-1$ valued matrix $M$, $O_{M}$ refers to the corresponding Cuntz-Krieger algebra. Their result
is a special case of a more general one, in which the stable
and unstable Ruelle algebras $R^{s}$ and $R^{u}$ associated to a
hyperbolic dynamical system are shown to be dual (see \cite{KPS}).

A particular example of a hyperbolic dynamical system is provided by
an Anosov diffeomorphism of a compact manifold; thus the duality
discovered by Kaminker and Putnam holds for these. An obvious question
is whether or not the same duality holds for Anosov $flows$. The
principal example of such a flow is given by 
geodesic flow on a compact, negatively curved
Riemannian manifold $M$. The algebras $R^{s}$ and $R^{u}$ can in this
case be regarded as foliation algebras. Specifically, define two
equivalence relations on $SM$ by respectively $v \sim_{s}w$ if $\limsup_{t \to \infty}
d_{SM}(g_{t}v,g_{t}w) = 0$, and $v \sim_{u} w$ if $\limsup_{t \to
  - \infty}
d_{SM}(g_{t}v,g_{t}w) = 0$. Define weak versions of these equivalence
relations by respectively $v \sim_{ws}w$ if $g_{t}(v) \sim_{s} w$ for
some $t$, and similarly for $v \sim_{wu} w$. The equivalence classes
of these latter two relations make up two codimension-$1$ foliations $\mathcal{F}^{ws}$ and
$\mathcal{F}^{wu}$ of $SM$. We can then form (see e.g. \cite{Co1}) the
corresponding foliation algebras $C^{*}_{r}(\mathcal{F}^{ws})$ and $C^{*}_{r}(\mathcal{F}^{wu})$. The work of
Kaminker and Putnam then suggested that $C^{*}_{r}(\mathcal{F}^{ws})$
should be dual in the aforementioned sense to $C^{*}_{r}(\mathcal{F}^{wu})$.

Now it is easy to see that the unit tangent sphere at a point of $M$
acts as a transversal to both foliations. We may therefore reduce the two
holonomy groupoids to this transversal and so obtain equivalent
groupoids, which are now $r$-discrete. Finally, it is easy to see that
these groupoids are in fact the same, and can be each identified with
the transformation groupoid $\partial \tilde{M} \rtimes \Gamma$,
where $\Gamma = \pi_{1}(M)$ and the boundary $\partial \tilde{M}$ is that
associated to the Gromov hyperbolic metric space $\tilde{M}$, acted
apon by $\Gamma$ by an extension of the action of $\Gamma$ by deck
transformations on $\tilde{M}$. Since $M$ is
compact and negatively curved, the group $\Gamma$ is of course hyperbolic in the sense of
Gromov, and $\partial\tilde{M}$ can be equivariantly identified with
$\partial \Gamma$. Consequently, if $C_{r}^{*}(\mathcal{F}^{ws})$ is to
be dual to $C_{r}^{*}(\mathcal{F}^{ws})$, we expect that the strongly
Morita equivalent algebra $C(\bgamma)
\rtimes \Gamma$ will be then dual to itself, or, equivalently, to its
opposite algebra. In other words, we can
reformulate the question of duality for the foliation algebras purely geometric-group-theoretically as follows:  is $C(\bgamma) \rtimes
\Gamma$ a Poincar\'e duality algebra when $\Gamma = \pi_{1}(M)$, for a
compact, negatively curved manifold $M$?

It is not difficult to see
that the answer to this question is yes in the case where $M$ has constant negative curvature. For then, if
say $n=2$ for simplicity, we may take $\Gamma$ to be a uniform lattice
in $G = PSL_{2}(\R)$, and then for $P$ equal to the parabolic subgroup of
upper triangular matrices of determinant $1$, we may identify $SM$
with $G/\Gamma$ and $\bgamma$ with $G/P$. Since the groupoids $G/P
\rtimes \Gamma$ and $G/\Gamma \rtimes P$ are 
equivalent, and since by two applications of the Thom Isomorphism,
$C(G/\Gamma) \rtimes P$ is $KK$-equivalent to $C(G/\Gamma) \cong
C(SM)$,  we see $C(\bgamma) \rtimes \Gamma$ is $KK$-equivalent to
$C(SM)$. Since $SM$ is a spin$^{c}$ manifold, $C(SM)$ has Poincar\'e
duality in $K$-theory, and therefore so does $C(\bgamma) \rtimes
\Gamma$.

Similar arguments can be used for the higher dimensional
cases of constant negative curvature. On the other hand, if the
curvature is variable, it seems to be necessary to use the infinite dimensional
techniques of Higson, Kasparov and Tu (\cite{Tu}). One then argues as
follows. The Baum-Connes conjecture for
the amenable groupoid $\bgamma \rtimes \Gamma$
tells us that $C(\bgamma) \rtimes \Gamma$ is
$KK$-equivalent to $C_{0}(\bgamma \times \ebar) \rtimes \Gamma \cong
C_{0}(S\tilde{M}) \rtimes \Gamma$ which in turn is strongly Morita
equivalent to $C(SM)$. Again, as $SM$ is a compact spin$^{c}$ manifold,
$C(SM)$ has Poincar\'e duality, and we are done.

These arguments do not however provide a concrete description of the fundamental
class $\Delta$, which
is desirable at least from the point of view of noncommutative
geometry (whose basic data are cycles, not merely classes.) To find such a concrete description was in fact the starting point
of our investigation. We wished, moreover, to describe such a cycle, purely in terms of the action of $\Gamma$ on its Gromov
boundary and without reference to spin$^{c}$ manifolds, Dirac
operators, and so on. That such a description exists was suggested by
the following example, of quite a different type from the above. 

Let
$\Gamma = \F_{2}$. Then $\Gamma$ is a hyperbolic group, with boundary
a Cantor set. It is easy to check (see e.g. \cite{Spiel}) that $C(\bgamma) \rtimes \Gamma$ is
in fact isomorphic to a Cuntz-Krieger algebra $O_{M}$  with matrix $M$ symmetric. By
the results of Kaminker and Putnam, we conclude for
reasons having apparently nothing to do with topology (but instead with the combinatorics of subshifts of finite type) that $C(\bgamma) \rtimes \Gamma$ is a Poincar\'e
duality algebra. For in this case $O_{M} \cong O_{M^{t}}$. Similar calculations verify that $C(\bgamma) \rtimes
\Gamma$ is a Poincar\'e duality algebra when $\Gamma$ is a free product of cyclic
groups.

Motivated by the latter calculations, we will in this paper approach
the problem from a different point of view, which will turn out to be
quite fruitful, yielding a Poincar\'e duality result for a very wide class of
hyperbolic groups, where neither the argument above in the case of $\Gamma =
\pi_{1}(M)$ nor that of Kaminker and Putnam appear (as far as
we know) to apply.

Let then $\Gamma$ be an arbitrary hyperbolic group and $A = C(\bgamma) \rtimes \Gamma$ the
corresponding cross product. Our method is as follows. We
construct a canonical extension of $A \otimes A^{\rm op}$ by the compact
operators based on
simple considerations of the action of the group $\Gamma$ on its compactification
$\bar{\Gamma}$. Specifically, associated to the
compactification, there are two extensions of the algebra $C(\bgamma)
\rtimes \Gamma$ by the compact operators, one corresponding, roughly,
to the action of $\Gamma$ on $l^{2}\Gamma$ by left translation and the
action of $C(\bar{\Gamma})$ by multiplication operators, and the
other to the action of $\Gamma$ by right translation and the action of
$C(\bar{\Gamma})$ by multiplication operators twisted by inversion on the group. Each extension yields a map $C(\bgamma)
\rtimes \Gamma$ into the Calkin algebra, and these two maps into the
Calkin algebra commute as
a consequence of the compactification being, in the language of
\cite{Hig2}, `good,' which simply means that metric balls of uniform size
become small in the topology of the compactification near the
boundary. Using this asymptotic commutativity, we obtain a single map from
$A \otimes A^{\rm op}$ into the Calkin algebra; i.e. an extension of
$A \otimes A^{\rm op}$ by the compact operators. We
define $\Delta$ to be the corresponding $KK$-class.

We will then set about proving that the class $\Delta \in KK^{1}(A \otimes
A^{\rm op} , \C)$ induces Poincar\'e duality, provided $\Gamma$ is
torsion-free and a certain condition regarding geodesics is
met. The latter can be stated as: the boundary has a continuous self
map with no fixed points; it is needed for a selection argument in
the latter stages of the proof. This 
technical condition is of course satisfied by groups whose boundaries are spheres
or Cantor sets; it is unknown to the author
whether there are any groups whose boundaries do not satisfy
it.  In our argument we will still make use of the Baum-Connes conjecture for
the groupoid $\bgamma \rtimes \Gamma$, but this time not to produce a
class which $a$ $priori$ we know induces Poincar\'e duality, as in the
discussion of $\Gamma = \pi_{1}(M)$ above, but to
show that our class $\Delta$ does.

The first step in proving that product with $\Delta$ does indeed
induce a Poincar\'e duality isomorphism, is to construct an
inverse, or dual element $\Dudelta \in KK^{-1}( \C , A \otimes A^{\rm
  op})$. We do this using a construction of Gromov, which produces a
sort of geodesic flow for an arbitrary hyperbolic group. We then show that $\Dudelta \otimes_{A^{\rm op}}\Delta =
1_{A}$. A calculation in \cite{Heath}
showed that in the case of the free group $\F_{2}$, the cycle
corresponding to the product
$\Dudelta \otimes_{A^{\rm op}} \Delta$ was a compact perturbation of
the ``$\gamma$-element'' cycle constructed by Julg and Vallette in \cite{JV},
parameterised by the points of $\bgamma$. In other words in this case the
statement $\Dudelta \otimes_{A^{\rm op}} \Delta = 1$ was equivalent to
the statement $\gamma_{\partial \F_{2} \rtimes \F_{2}} = 1_{C(\partial
  \F_{2}) \rtimes \F_{2}}$ where $\gamma_{\partial \F_{2} \rtimes
  \F_{2}}$ is the $\gamma$-element for this transformation groupoid,
and so roughly equivalent to the statement that the Baum-Connes map for the groupoid is
an isomorphism. The latter has been verified by Tu (\cite{Tu}) for general hyperbolic
groups, and we are able to resolve the general case in a somewhat
analogous way.

The organization of the paper is as follows. In Section 2 we provide a
summary of the basic facts from $KK$-theory which we will need. In
Section 3 we set up the formalism of $K$-theoretic Poincar\'e
duality. In Section 4 we construct the fundamental class $\Delta$,
which as mentioned
exists
for every hyperbolic group, with or without torsion, and with or
without a fixed-point-free map on the boundary. We then
construct the dual element $\Dudelta$ using an analog for hyperbolic
groups of geodesic flow on a negatively curved manifold. In Section
5 we begin the process of verifying the fundamental equation of
Poincar\'e duality: $\Dudelta \otimes_{A^{\rm op}}\Delta = 1_{A}$,
where $A = C(\bgamma) \rtimes \Gamma$.

Given the class $\gamma_{A} = \Dudelta
\otimes_{A^{\rm op}}\Delta \in KK(A,A)$, we wish to show it is
$1_{A}$. We first calculate the cycle corresponding to $\gamma_{A}$. We then make use of this calculation to show that $\gamma_{A}$ lies in the range of the descent map
$$\lambda: RKK_{\Gamma}(\bgamma ; \C , \C) \to KK(A,A).$$We 
reduced to showing that its preimage, $\gamma_{\bgamma}$, is
$1_{\bgamma} \in RKK(\bgamma ; \C , \C)$, since the descent map has
the property that $\lambda (1_{\bgamma}) = 1_{A}$. The Baum-Connes
conjecture for the amenable groupoid $\bgamma \rtimes \Gamma$ implies that there is an isomorphism
$RKK_{\Gamma}(\bgamma ; \C , \C) \cong RKK_{\Gamma}(\bgamma \times
\ebar ; \C , \C)$, where $\ebar$ is the classifying space for proper
actions of $\Gamma$,  and so it suffices to show that the image of
$\gamma_{\bgamma}$ under this isomorphism is $1_{\bgamma \times
  \ebar}$. This calculation, which though not difficult is slightly
involved, is performed in Sections 6 and 7. It is at this
point that we require the hypothesis that the boundary of $\Gamma$
possesses a fixed-point-free map.

I would like to thank N. Higson, my advisor from the Pennsylvannia
State University, as well as J. Kaminker and I. Putnam,
for extremely valuable comments and suggestions regarding the material
in this paper. Finally, I would like to thank the referees, for several useful
remarks.

\section{KK-theoretic preliminaries}

Kasparov's $KK$-theory, along with some of its elaborations, will be
used extensively in this paper. $KK$ can be understood
categorically (\cite{Hig4}). From this latter point of view, there is a category
$\mathbf{KK}$ whose objects are separable, nuclear $C^{*}$-algebras and whose morphisms
$A \to B$ are the elements of $KK(A,B)$. There is a functor from
the category of $C^{*}$-algebras to the category $\mathbf{KK}$. There
is a composition, or intersection product operation $KK(A , D)
\times KK(D , B) \to KK(A,B)$ which we denote by $(\alpha , \beta)
\mapsto \alpha \otimes_{D} \beta$. If $D$ is a $C^{*}$-algebra, there is a natural map $KK(A,B) \to
KK(A \otimes D , B \otimes D)$, $\alpha \mapsto \alpha
\otimes 1_{D}$, and similarly a map $KK(A,B) \to KK(D \otimes A , D
\otimes B)$. The above three operations imply the existence of a mixed
cup-cap product $$KK(A _{1} , B_{1} \otimes D)
\times KK(D \otimes B_{2} , A_{2}) \rightarrow
KK(A_{1} \otimes B_{2} , B_{1} \otimes A_{2})$$which is denoted
$(\alpha , \beta) \mapsto \alpha \otimes_{D} \beta$, and
defined by $\alpha \otimes_{D}\beta = \bigl(\alpha \otimes 1_{B_{2}}\bigr) \;
\otimes_{B_{1} \otimes D \otimes B_{2}}\; \bigl(1_{B_{1}}
\otimes \beta \bigr)$. There are higher $KK$ groups $KK^{i}(A,B)$ for all $i \in
\Z$, defined by $KK^{i}(A , B) = KK(A , B \otimes C_{i})$ where
$C_{i}$ is the $i$th complex Clifford algebra, and one of the features of the
theory is that the intersection product is graded commutative. If $A_{1},  \ldots ,  A_{n}$ are
$C^{*}$-algebras, let $\sigma_{ij}$
denote the map $$A_{1}\otimes \cdots  A_{i} \otimes
\cdots A_{j} \otimes \cdots \otimes A_{n} \rightarrow
A_{1}\otimes \cdots  A_{j} \otimes \cdots A_{i}
\otimes \cdots \otimes A_{n}$$ obtained by flipping the
two factors. Then by graded commutativity we mean:

\begin{lemma}

If $\alpha \in KK^{i}(A_{1},B_{1})$ and $\beta \in KK^{j}(A_{2} ,
B_{2})$, then $$\alpha \otimes_{\C} \beta = (-1)^{ij} \,
(\sigma_{12})_{*} \sigma_{12}^{*}( \beta \otimes \alpha) \in KK(A_{1}
\otimes  A_{2} , B_{1} \otimes B_{2}).$$

\end{lemma}

Let $\Lambda$ be a discrete group.Then as well as the category
$\mathbf{KK}$ there is the category $\mathbf{KK_{\Lambda}}$, whose
objects are $\Lambda - C^{*}$-algebras and whose morphisms are the
elements of $KK_{\Lambda}(A,B)$. We can think of these as equivariant
morphisms. There is a descent map $$\lambda : KK_{\Lambda}(A , B) \to
KK(A \rtimes \Lambda , B \rtimes \Lambda)$$producing from an
equivariant morphism a nonequivariant one. There is a map backwards if
$A$ and $B$ happen both to be trivial $\Lambda - C^{*}$-algebras in the
sense that every $\gamma \in \Lambda$ acts as the identity
automorphism. The descent map is natural: that is, $\lambda ( \alpha
\otimes_{D} \beta)  = \lambda (\alpha) \otimes_{D \rtimes \Lambda}
\lambda (\beta)$. The group $KK_{\Lambda}(A,A)$ is a ring with the
intersection product, and there is an identity in this ring, denoted
$1_{A}$, and it satisfies $\lambda (1_{A}) = 1_{A \rtimes \Lambda}$.

Finally, let $X$ be a locally compact $\Lambda$ space. Then there is
another category, denoted $\mathbf{\mathcal{RKK}_{\Lambda}}$, this time
whose objects are $\Lambda - C(X)$-algebras (see \cite{Ka1}) and whose
morphisms are the elements of $\mathcal{RKK}_{\Lambda}(X ; A , B)$. In the case of $A
= C_{0}(X) \otimes A_{0}$ and $B = C_{0}(X) \otimes B_{0}$, with $A_{0}$ and
$B_{0}$ $\Lambda - C^{*}$-algebras, we denote, following Kasparov,
the group $\mathcal{RKK}_{\Lambda}(X; A , B)$ by $RKK_{\Lambda}(X ;
A_{0} , B_{0})$. The intersection product
$$\mathcal{RKK}_{\Lambda}(X; A , D) \times \mathcal{RKK}_{\Lambda}(X; D
, B) \to \mathcal{RKK}_{\Lambda}(X; A , B)$$is denoted $(\alpha ,
\beta) \mapsto \alpha \otimes_{X , D} \beta$, and similarly for
$RKK_{\Lambda}$. Note also that $RKK_{\Lambda}(X ; A , A)$ has
a unit, which is denoted $1_{X , A}$, and if $A = \C$ we denote this
unit simply by $1_{X}$. Finally, if $Z$ is any space, there is a natural map
$$p_{Z}^{*}: RKK_{\Lambda}(X; A , B) \to \mathcal{RKK}_{\Lambda}(X \times Z
  ;A , B).$$ This map is natural with respect to intersection
  products and thus is a ring homomorphism when $A = B$. Under certain special circumstances it is an
  isomorphism (see Theorem 54).

Throughout this paper we will let $B(\mathcal{E})$ denote bounded
operators on a Hilbert
module $\mathcal{E}$, $K(\mathcal{E})$ compact operators, and $Q(\mathcal{E})$ the Calkin algebra
$B(\mathcal{E}) / K(\mathcal{E})$. The projection map $B(\mathcal{E})
\to Q(\mathcal{E})$, which will be invoked frequently, will always be
denoted by $\pi$. 

Following Kasparov (\cite{Ka1}), if $\mathcal{E}$ is a Hilbert
$B$-module and $A$ acts on $\mathcal{E}$ by a homomorphism $A \to
B(\mathcal{E})$, we will refer to $\mathcal{E}$ as a Hilbert $(A,B)$-bimodule.

Because all the algebras in this paper are
ungraded -- or alternatively, have trivial grading --  we can make certain simplifications in the
definitions of the $KK$ groups (see \cite{Bla}). With such ungraded
$A$ and $B$, cycles for $KK(A,B)$ are given simply by pairs
$(\mathcal{E} , F)$ where $\mathcal{E}$ is an $(A,B)$-bimodule, $F$ commutes modulo compact operators with the action of
$A$, and $a(F^{*}F - 1)$ and $a(FF^{*}-1)$ are compact for every $a \in A$.

Cycles for $KK^{1}(A,B)$ are given  by pairs $(\mathcal{E} , P)$ for
which $P$ is as before an operator on the $(A,B)$-bimodule
$\mathcal{E}$ as above, and where $P$ satisfies the three conditions
$[a,P]$, $a(P^{2} - P)$, and $a(P-P^{*})$
are compact for all $a \in A$. Such pairs are equivalently given by
$extensions$, i.e. homomorphisms $A \mapsto Q(\mathcal{E})$. For by
the Stinespring construction, under our nuclearity assumptions, for
each such homomorphism $\tau$ there exists a Hilbert $(A,B)$-module
$\tilde{\mathcal{E}}$, an isometry $U: \mathcal{E} \to
\tilde{\mathcal{E}}$, and an operator $P$ on $\tilde{\mathcal{E}}$
such that $a(P^{2} - P)$, $[a,P]$, and $a(P-P^{*})$ are compact for all $a \in A$,
and $\pi (U^{*}PaPU) = \tau (a)$ for all $a \in A$.

Recall that $KK^{-1}(\C , C^{*}(\R)) \cong \Z$ and is generated
by the class $[\hat{d}_{\R}]$ of the Dirac operator on $\R$, viewed as an unbounded self-adjoint
multiplier of $C^{*}(\R)$. The class
$[\hat{d}_{\R}]$ allows us to identify, for any $C^{*}$-algebras $A$
and $B$, the groups $KK^{1}(C^{*}(\R) \otimes A , B)$, and
$KK(A,B)$, by the map $KK^{1}(C^{*}(\R) \otimes A , B) \to KK(A,B)$,
$x \mapsto [\hat{d}_{\R}]\otimes_{C^{*}(\R)}
x$. We shall need to compute this map at the level of cycles in
several simple cases.

 Let $\psi$ be the function in $C^{*}(\R)$ whose Fourier
transform is $\frac{-2i}{z+i}$. It has the property that $ \psi + 1$
is unitary in $C^{*}(\R)^{+}$.

\begin{lemma}

Let $A$ be a $C^{*}$-algebra, $\varphi$ a homomorphism $C^{*}(\R) \to
A$, and suppose $\tau: A \to Q(H)$ is a homomorphism to the Calkin
algebra. Let
$[\tau]$ denote the class in $KK^{1}(A , \C)$ corresponding to
$\tau$. Then the class $[\hat{d}_{\R}] \otimes_{C^{*}(\R)} \varphi^{*}([\tau])
\in KK(\C , \C)$ is represented by the
cycle $(H , T+1)$, where $T$ is any operator on $H$ such that $\pi (T)
= \tau (\varphi (\psi))$. 

\end{lemma}

We will also need the following simple lemma.

\begin{cor}

Define a class $[\tau] \in KK^{1}(C^{*}(\R) , \C)$ by means of the homomorphism $\tau: C^{*}(\R)
\rightarrow Q(L^{2}(\R))$, $$f \mapsto \pi \bigl(\chi \cdot \lambda
(f)\bigr),$$ where $\lambda$ is the left regular representation of
$C^{*}(\R)$ and $\chi$ is the characteristic function of the left
half-line. Then $[\hat{d}_{\R}] \otimes_{C^{*}(\R)}
 [\tau] =
 [ 1_{\C}] \in KK(\C, \C)$.

\end{cor}

\begin{proof}

This follows from Lemma 2 and a calculation; one checks simply that $\chi \cdot \psi$ as an operator on $L^{2}(\R)$ has index 1. One can do this by solving a simple differential equation. (See \cite{Ska}).

\end{proof}

\begin{note}

\rm

Remark that the function $\chi$ above can be replaced by any function
on $\R$ which is $1$ at $-\infty$ and $0$ at $+ \infty$. For any such
function gives the same extension.

\end{note}

Next, let $A_{1}$ and $A_{2}$ be $\Lambda - C^{*}$-algebras, where $\Lambda$ is a
discrete group. An action of $\Lambda$ on an
$(A_{1},A_{2})$-bimodule $\mathcal{E}$ will always refer to an action of
$\Lambda$ as complex linear maps compatible with the inner product in
the sense that $<\gamma \xi , \gamma \eta >_{A_{2}} = \gamma ( <\xi ,
\eta>_{A_{2}})$, and compatible with the bimodule structure in the sense
that $\gamma (a\xi b) = \gamma (a) \gamma (\xi) \gamma (b)$. Such
$\mathcal{E}$ will be referred to as a $\Lambda -
(A_{1},A_{2})$-bimodule. If we wish to possibly waive the part of the
last requirement that states that $\gamma (a \xi) = \gamma (a) \gamma
(\xi)$, whilst maintaining the requirement that $\gamma (\xi b) = \gamma
(\xi) \gamma(b)$, we will simply call $\mathcal{E}$ a
$\Lambda-A_{2}$-module. Thus, such a module satisfies $\gamma ( \xi b)
= \gamma (\xi) \gamma (b)$, but the homomorphism $A_{1} \to
\mathcal{B}(\mathcal{E})$ may not necessarily be
$\Lambda$-equivariant

Cycles for $KK_{\Lambda}(A_{1},A_{2})$ are then given
by pairs $(\mathcal{E} , F)$ where $\mathcal{E}$ is a $\Lambda -
(A_{1},A_{2})$-bimodule, and where $F \in B(\mathcal{E})$ with $a(F^{*}F - 1)$
and $a(FF^{*} - 1)$ compact for all $a \in A_{1}$, and $\gamma (F) - F$
compact for all $\gamma \in \Lambda$. Cycles for $KK_{\Lambda}^{1}(A_{1} ,
A_{2})$ are given by pairs $(\mathcal{E} , P)$ where $\mathcal{E}$ is a
$\Lambda - (A_{1},A_{2})$-bimodule and $P$ is an operator with
$a(P^{2} - P)$, $a(P-P^{*}),$ and $[a,P]$ compact for all $a \in A_{1}$, and $\gamma (P)
- P$ compact for all $\gamma \in \Lambda$. 

A minor technical issue which in general we do not know how to resolve
concerns the question of whether or not an equivariant map $A_{1} \to
Q(\mathcal{E})$, where $\mathcal{E}$ is a $\Lambda - A_{2}$-module,
produces an element of $KK^{1}_{\Lambda}(A_{1} , A_{2})$. If $\Lambda$ is the
trivial group this is of course the Stinespring
construction, given our standing assumption that all $C^{*}$-algebras
(with the obvious exceptions of Calkin algebras and so on) are nuclear. In the general case, an equivariant homomorphism $A_{1} \to Q(\mathcal{E})$ yields a homomorphism $A_{1} \rtimes
\Lambda \to Q(\mathcal{E} \rtimes \Lambda)$ where $\mathcal{E} \rtimes
\Lambda$ is as in \cite{Ka1}, being a certain $(A_{1}\rtimes \Lambda , A_{2}
\rtimes \Lambda)$-bimodule (this is part of the definition of the
descent map) and so an element of $KK^{1}(A_{1}
\rtimes \Lambda , A_{2} \rtimes \Lambda)$ as long as not merely
$A_{1}$ and $A_{2}$ are nuclear, but also $A_{1} \rtimes \Lambda$ and
$A_{2} \rtimes \Lambda$ are nuclear. But such an element may not
necessarily come under descent from an element of $KK^{1}_{\Lambda}
(A_{1}, A_{2})$. To avoid this issue, we make the following definition.

\begin{defn} \rm

Let $\Lambda$ be a discrete group, let $A_{1}$ and $A_{2}$ be $\Lambda -
C^{*}$-algebras and let $ \mathcal{E}$ be a $\Lambda - A_{2}$-module. Let $\tau: A_{1} \to Q(\mathcal{E})$ be a
$\Lambda$-equivariant homomorphism. We say $\tau$ is dilatable if
there is a $\Lambda - (A_{1} , A_{2})$-bimodule $\tilde{\mathcal{E}}$,
an operator $P$ on $\tilde{\mathcal{E}}$ such that
$[a,P]$, $a (P^{2} - P)$, $a(P^{*} - P)$ and $\gamma (P) - P$ are compact for all $a
\in A_{1}$, $\gamma \in \Lambda$, and if there exists an isometry
$U: \mathcal{E} \to \tilde{\mathcal{E}}$, such that $\pi (U^{*}PaPU) = \tau (a) \in Q(\mathcal{E})$ for all $a \in
A_{1}$.

\end{defn}

As mentioned above, if $\Lambda$ is the trivial group then every
homomorphism $A_{1} \to Q(\mathcal{E})$ is dilatable. The same is clearly
true of finite $\Lambda$. In general, with the hypothesis of
dilatibility, we do clearly have the following: 

\begin{lemma}

If $A_{1}$, $A_{2}$, $\mathcal{E}$, $\Lambda$ and $\tau$ as
above, and if $\tau$ is dilatable, then $\tau$ defines 
a class $[\tau]$ in $KK^{1}_{\Lambda}(A_{1} , A_{2})$ by the pair
$(\tilde{\mathcal{E}} , P)$.

\end{lemma}

We next pass to a case where to calculate the Kasparov product of two
elements one of which is given by a dilatable homomorphism, we do not
need to explicitly involve the dilation. We will use this technical lemma
several times, sometimes with $\Lambda$ the trivial group. In the
latter case, the lemma gives a method of avoiding explicit
construction of a completely positive section.

\begin{lemma}

 Let $A_{1} , A_{2}$ be $\Lambda - C^{*}$-algebras and $\mathcal{E}$
 be a $\Lambda - A_{2}$-module. Let $[h]$ be a class in $KK^{1}_{\Lambda}(C^{*}(\R)
\otimes A_{1} , A_{2})$ given by 
a $\Lambda$-equivariant dilatable homomorphism $h:  C^{*}(\R) \otimes A_{1} \to Q(\mathcal{E})$ of
the form $x \otimes a_{1} \mapsto
h'(x) h''(a_{1})$, where $h'$ and $h''$ are $\Lambda$-equivariant
homomorphisms. Suppose that the homomorphism $h''$ lifts
to a $\Lambda$-equivariant homomorphism $\tilde{h}'': A_{1} \to B(\mathcal{E})$. Then
the class $[\hat{d}_{\R}] \otimes_{C^{*}(\R)} [h] \in KK_{\Lambda}(A_{1},A_{2})$ is
represented by the
following cycle. The module is $\mathcal{E}$ with its
 original $\Lambda - A_{2}$-module structure and the left $A_{1}$-module
structure given by the homomorphism $\tilde{h}''$. The operator is
 given by $F+1$ where $F$ is any operator on
 $\mathcal{E}$ such that $\pi (F) =
 h'(\psi)$.

\end{lemma}

\begin{rmk}

\rm

Similar lemmas can be formulated and proved for the $RKK_{\Lambda}$
category, but we leave it to the reader to formulate them. 

\end{rmk}

\section{Formalism of Noncommutative Poincar\'e Duality}

Let us begin with a lemma. See \cite{KP} for a similar discussion.

\begin{lemma}

Let $A$ and $B$ be $C^{*}$-algebras and let $\Delta$ and $\Dudelta$ be
two elements in $KK^{i}( A\otimes B , \C)$ and $KK^{-i}(\C , A\otimes
B)$ respectively. Define a map $\Dudelta_{j}: K^{j}(B) \mapsto
K_{j-i}(A)$ by $\Dudelta_{j}(x)= \Dudelta \otimes_{B} x $. Define a map
$\Delta_{j}: K_{j}(A) \mapsto K^{j +i}(B)$ by $\Delta_{j}(y) =
y\otimes_{A} \Delta.$ Define also two classes in
respectively $KK(A,A)$ and $KK(B,B)$ by $\gamma_{A} = \bigl( \Dudelta \otimes 1_{A} \bigr) \otimes_{A \otimes
  B \otimes A} \bigl( 1_{A} \otimes \sigma_{12}^{*}(\Delta) \bigr)$,
and $\gamma_{B} = \bigl( (\sigma_{12})_{*}(\Dudelta ) \otimes 1_{B}
\bigr) \otimes_{B \otimes A \otimes B} \bigl( 1_{B} \otimes \Delta
\bigr).$ Then we have:
$$\Delta _{j-i}( \Dudelta _{j} (x)) = (-1)^{ij}\,\gamma_{B}
\otimes_{B} x , \; x \in K^{j}(B);$$ and
$$\Dudelta_{j+i}(\Delta_{j}(y)) = (-1)^{ij}\, y \otimes_{A}
\gamma_{A}, \; y \in K_{j}(A).$$

\end{lemma}

\begin{proof}
We verify the first equation; the second follows similarly. Let $x \in K^{j}(B).$ Then it follows from the definition that
$$\Delta_{j-i}(    \Dudelta_{j}(x))= (\Dudelta \otimes 1_{B} )
\otimes_{A\otimes B\otimes B} ( 1_{A}\otimes   x \otimes 1_{B})
\otimes_{A\otimes B} \Delta.$$By functoriality of the intersection
product we may write this $$ \bigl((\sigma_{12})_{*}(\Dudelta ) \otimes 1_{B} \bigr)
\otimes_{ B\otimes A \otimes B} \sigma_{12}^{*}( 1_{A}\otimes   x \otimes 1_{B})
\otimes_{A\otimes B} \Delta.  $$On the other hand, again by
definition, we have $$\gamma_{B}\otimes_{B} x = \bigl((\sigma_{12})_{*}(\Dudelta) \otimes
1_{B}\bigr) \otimes_{B\otimes A \otimes B} (1_{B} \otimes \Delta)
\otimes_{B} x.$$ So we are reduced to proving that $(1_{A} \otimes x \otimes
1_{B}) \otimes_{A \otimes B}\Delta = (-1)^{ij}(1_{B} \otimes \Delta) \otimes_{B}
x$. But this follows immediately from Lemma 1.

\end{proof}

In view of this theorem, we will take as the definition of duality between two $C^{*}$-algebras the following (compare \cite[page 588]{Co2}):

\begin{defn}\rm
Two separable, unital, and nuclear $C^{*}$-algebras $A$ and $B$ are
$dual$ with a dimension shift of $i$ if there exists $\Delta \in
KK^{i}(A\otimes B, \mathbb C)$, $\Dudelta \in KK^{-i}(\C , A  \otimes
B )$ such that $$\Dudelta \otimes_{B}\Delta = 1_{A}$$ and $$\Dudelta
\otimes_{A} \Delta = (-1)^{i}\, 1_{B}.$$ We
will call such a pair $(\Dudelta , \Delta)$ a duality pair.

\end{defn}

\begin{thm}

If $A$ and $B$ are dual in the sense of Definition 10, then the maps
$\Dudelta_{*}$ and $\Delta_{*}$ defined in Lemma 9 induce inverse isomorphisms up to the signs specified there $K_{j}(A) \cong K^{j+i}(B)$ and $K^{j}(B) \cong K_{j-i}(A)$.

\end{thm}

 For the next piece of terminology recall that for a $C^{*}$-algebra
$A$, $A^{\rm op}$ denotes the opposite algebra of $A$. 

\begin{defn}
\rm

A separable, nuclear $C^{*}$-algebra $A$ is a Poincar\'e duality algebra
 if $A$ and $A^{\rm op}$ are dual in the sense of
Definition 10. We will refer to $\Delta$ as the Fundamental class of
$A$, and $(\Dudelta , \Delta)$ as a Poincar\'e duality pair.

\end{defn}

\section{The Main Theorem}

Let $\Gamma$ be a hyperbolic group. We shall assume here and
throughout this paper that $\Gamma$ is torsion-free. To $\Gamma$ we can add a
boundary $\bgamma$ which compactifies the group $\Gamma$ understood as
a metric space. Thus, $\bar{\Gamma} = \Gamma \cup
\bgamma$ can be given the structure of a compact metrizable space in
which $\Gamma$ sits densely. For details see \cite{GH}. The group $\Gamma$ acts by homeomorphisms on $\bgamma$ and
this action is topologically amenable in the sense of \cite{Del} (see
the appendix of \cite{Del} for a proof of this fact). Therefore, to each
hyperbolic group we can associate an amenable $r$-discrete amenable groupoid $\bgamma \rtimes
\Gamma$ and then a groupoid $C^{*}$-algebra
$C(\bgamma) \rtimes \Gamma$ which for the rest of this paper we shall
denote  by $A$. The $C^{*}$-algebra $A$ is separable, simple, nuclear
and purely infinite (see \cite{Spiel} or \cite{Del2}). Our goal is
to show that
for a large subclass of hyperbolic groups $\Gamma$, $A$ is a
Poincar\'e duality algebra in the sense of Definition 12. Let us first
state certain simple facts we shall require.

\begin{note}\rm

When we are thinking of elements of $\Gamma$ as simply points in the
metric space $\Gamma$, we shall use the notation $x,
y$, etc. In particular, $x_{0}$ will always refer to the identity of
the group, viewed as a natural basepoint. Also, for
any $R \ge 0$ and any $x \in \Gamma$, $B_{R}(x)$ denotes the ball
of radius $R$ (with respect to the word metric) centered at $x$.

\end{note}

For
convenience we will also fix a metric $d_{\bar{\Gamma}}$ on
  $\bar{\Gamma}$ compatible with the topology. The
following lemma then follows from the definition of the
topology on $\bar{\Gamma}$ (see \cite{GH}).

\begin{lemma} 
 If $\epsilon >0,$ there exists $R \ge 0$ such that if $a,b \in
 \bar{\Gamma}$ and $d_{\bar{\Gamma}}(a,b) \ge \epsilon$, then every
 geodesic from $a$ to $b$ passes through $B_{R}(x_{0})$. Conversely,
 if $R \ge 0$, there exists $\epsilon >0$ such that if every geodesic
 between $a$ and $b$ passes through $B_{R}(x_{0})$, then $d_{\bar{\Gamma}}(a,b) \ge \epsilon$.

\end{lemma}

We will also require the following. Recall that $(x \; | \; y)$
denotes the Gromov product of $x,y \in \Gamma$ (see \cite{Gro} or
\cite{GH}). For the proof of this lemma see for example \cite{Roe}.

\begin{lemma}

If $f$ is a bounded function on $\Gamma$, then $f$ extends to a
continuous function on $\bar{\Gamma}$ if and only if for all $\epsilon
> 0$ there exists $R \ge 0$ such that if $(x \; | \; y ) > R$, then
$|f(x) - f(y) | < \epsilon$. 

\end{lemma}

We shall need an explicit description of the classifying space for
proper actions of $\Gamma$. This is given by the Rips construction.

\begin{defn}

\rm

The Rips complex for $\Gamma$ of parameter $N$, $P_{N}(\Gamma)$, is
the simplicial complex whose vertices are the points of $\Gamma$, and
whose $k$-simplices are the sets of cardinality $k$ of diameter less
than or equal to $N$.

\end{defn}

Let $\bar{P}_{N}(\Gamma)$ denote the realization of the Rips
Complex. It can be viewed as the collection of finitely supported
probability measures on $\Gamma$ whose support has diameter $\le
N$. This point of view will be useful later on the proof when some
linear interpolation will be needed from $\Gamma$ to
$\bar{P}_{N}(\Gamma)$. Note that $\Gamma$ is embedded naturally in $\bar{P}_{N}(\Gamma)$. Clearly $\bar{P}_{N}(\Gamma)$ carries a free, simplicial, isometric, proper, co-compact action of $\Gamma$.

A proof of the following may be found in \cite{Meintriup}.

\begin{lemma}

For large enough $N$, $\bar{P}_{N}(\Gamma)$ is the classifying space $\ebar$ for proper actions of $\Gamma$.

\end{lemma}

\begin{note} \rm From this point onwards, we fix $N$ sufficiently large as
 in the above lemma, and denote the realization of the Rips complex
 with parameter $N$ simply by $\ebar$. We will also fix a simplicial
 metric $d_{\ebar}$ on $\ebar$, so that $\Gamma$ is isometrically embedded in
 $\ebar$ as the vertices of the complex. Then $\ebar$ is of course a
 hyperbolic space in its own right, and is quasi-isometric to
 $\Gamma$. 

\end{note}

We now pass to the construction of the fundamental class $\Delta \in
 KK^{1}(A \otimes A^{\rm op} , \C)$, which will arise naturally as an
 $extension$, or equivalently as a homomorphism $A \otimes A^{\rm op}
 \to Q(H)$ for some Hilbert space $H$. This map $A \otimes A^{\rm op}
 \to Q(H)$ will be given by two commuting maps $A \to Q(H)$ and $A^{\rm op}
 \to Q(H)$, which we shall denote by $\lambda$ and $\lambda ^{\rm op}$
 respectively.

 Passing to the description of $\lambda$, let us put $H =
 l^{2}(\Gamma)$. This notation will be retained throughout the rest of
 this paper. Let $e_{x}$, $x \in \Gamma$  denote the standard
 basis element of $H$ corresponding to point mass at $x$. For $\gamma \in \Gamma$ let $u_{\gamma}$
 denote the unitary in $B(H)$ given by left translation by $\gamma$, i.e.
 $u_{\gamma}( e_{x}) = e_{\gamma x}$. Let $\lambda (\gamma)$ denote
 the image of $u_{\gamma}$ in the Calkin algebra. Let $f$ be a
 function in $C(\bgamma)$, apply the Tietze extension theorem to
 extend $f$ to a continuous function $\tilde{f}$ on
 $\bar{\Gamma}$, and let $\lambda (f)$ denote the image in $Q(H)$ of
 the operator on $H$ given by multiplication by
 $\tilde{f}$, in other words the operator
 $e_{x} \mapsto  \tilde{f}(x)e_{x}$. Remark that though the map
 $\gamma \to u_{\gamma}$, $f \to \tilde{f}$ is not well-defined into
 $B(H)$, it $is$ well-defined into $Q(H)$, since any two extensions of a
 function $f$ differ by a function vanishing at $\infty$ and thus by a
 compact operator. The following lemma is a trivial calculation:

\begin{lemma}

 The assignment $\gamma \mapsto \lambda (\gamma)$, $f \mapsto \lambda
 (f)$, defines a covariant pair for the $C^{*}$-dynamical system
 $(C(\bgamma) , \Gamma)$, and so a homomorphism $A \to Q(H)$.

\end{lemma}

Next, define a map $\lambda ^{\rm op}: A^{\rm op} \to Q(H)$ as
follows. First, let $v_{\gamma}$, for $\gamma \in \Gamma$, denote the
unitary operator of $right$ translation by $\gamma$: $v_{\gamma}
(e_{x}) = e_{x\gamma}$. Let $\lambda ^{\rm op}(\gamma)$ denote the
image of this unitary operator in the Calkin algebra. If now $f \in
C(\bgamma)$, let $\tilde{f}$ denote an extension of $f$ to a
continuous function on $\bar{\Gamma}$ as
before, and let $\lambda ^{\rm op}(f)$ denote the image in the Calkin
algebra of the multiplication operator given by multiplication by the
function $x \mapsto \tilde{f}(x^{-1})$. These two maps are similarly
well-defined into the Calkin algebra, and we have easily:

\begin{lemma}
 The assignment $\gamma \mapsto \lambda ^{\rm op}(\gamma)$, $f \mapsto \lambda ^{\rm op}(f)$, defines a covariant pair with respect to the opposite action of $\Gamma$ on $C(\bgamma)$ and hence a homomorphism $A^{\rm op} \to Q(H)$.

\end{lemma}

We next show the two homomorphisms $\lambda$ and $\lambda ^{\rm op}$
commute as maps into the Calkin algebra. This follows from the
following

\begin{lemma}
Let $\tilde{f}$ be a function on $\Gamma$, viewed as a multiplication operator on $H$, and let $\gamma \in \Gamma$. \\
(1) If $x \mapsto \tilde{f}(x)$ is continuous on $\bar{\Gamma}$, then $[v_{\gamma} , \tilde{f}]$ is a compact operator. \\
(2) If $x \mapsto \tilde{f}(x^{-1})$ is continuous on $\bar{\Gamma}$, then $[u_{\gamma} , \tilde{f}]$ is a compact operator. 

\end{lemma}

\begin{proof}

 Let $\tilde{f}$ be as in (1). Choose $\epsilon > 0 $. Remark if $x,
 \gamma \in \Gamma$ we have $(x,x\gamma) \ge |x| - |\gamma|$. From
 this and Lemma 15 we see: there exists $K$
 such that $|x| > K \Rightarrow |\tilde{f}(x) - \tilde{f}(x\gamma)| <
 \epsilon$. In other words, the function $\tilde{f}(x) -
 \tilde{f}(x\gamma)$ vanishes at infinity. It follows immediately that
 $v_{\gamma}\tilde{f}v_{\gamma ^{-1}} - \tilde{f}$ is compact; for
 this operator is precisely multiplication by this function. Hence $\bigl( v_{\gamma}\tilde{f}v_{\gamma ^{-1}} -
 \tilde{f}\bigr)v_{\gamma} = [v_{\gamma} , \tilde{f}]$ is also a
 compact operator. (2) follows from (1) by conjugating by the
 unitary $H \to H$ induced from inversion on the
 group.
\end{proof}

\begin{rmk}

\rm

The above lemma can be restated in a slightly more general 
way. Having fixed a left-invariant metric on $\Gamma$, as we have
done, $right$ translation by a fixed $\gamma \in \Gamma$ gives an
operator of finite propagation; on the other hand any operator
of finite propagation commutes modulo compacts with multiplication by
a function in $C(\bar{\Gamma})$ by the same proof as that of Lemma
21.

\end{rmk}

\begin{defn}

\rm

Let $\Gamma$ be any hyperbolic group and $\bgamma$ its Gromov
boundary. Let $H$ denote $l^{2}(\Gamma)$. We define the fundamental class of the $C^{*}$-algebra $A =
C(\bgamma) \rtimes \Gamma$ to be the class $\Delta$ in $KK^{1}(A \otimes
A^{\rm op} , \C)$ corresponding to the homomorphism $A \otimes A^{\rm op}
\to Q(H)$ induced by the two commuting homomorphisms
$\lambda$ and $\lambda ^{\rm op}.$

\end{defn}

\begin{rmk}

\rm

Let $\Gamma$ be a discrete, not necessarily hyperbolic group acting co-compactly and properly on a nonpositively curved
space $X$, and let $\partial X$ denote the visibility boundary of
$X.$ The visibility boundary compactifies the group $\Gamma$ and all
of the above constructions extend to this situation. We thus obtain a map $C(\partial X) \rtimes \Gamma \otimes_{\rm max} \bigl(C(\partial X)
\rtimes \Gamma \bigr)^{\rm op} \to Q(H)$ in the same way. However, as
the $\Gamma$-action on $\partial X$ is no longer amenable, it is no
longer necessarily the case that such a map defines a $KK^{1}$ element. 

\end{rmk}

\begin{rmk}

\rm

If $J$ denotes the conjugate linear operator $H \to H$ sending the
element
$\sum_{\gamma} a_{\gamma} e_{\gamma} \in C_{c}(\Gamma)$ to the element
$\sum_{\gamma}
\bar{a}_{\gamma} e_{\gamma ^{-1}}$, then the
equation $J\lambda (a^{*}) J^{-1} = \lambda^{\rm op} (a)$ holds for
any $a \in A$. This is the content of Connes' Reality axiom
(see \cite{Co1}), except that the relation holds in the Calkin algebra rather
than in $B(H)$. In fact, it is easy to see that all our constructions
are compatible with the various real structures on the algebras,
Hilbert spaces, and so on, concerned, and that the cycle $\Delta$
in actually gives a $KR$-homology class. Similarly we shall see that
$\Dudelta$ gives a $KR$ class, and that the duality we are going to prove
holds in the real as well as the complex setting. 

\end{rmk}

We now proceed to the element $\Dudelta$, to construct which we shall use an idea of
Gromov and subsequent work by Champetier
and Matheus. Theorem 27 was first stated by Gromov (see \cite{Gro}, pg. 222), with a sketch of a
proof; details were added by the latter two authors in respectively
\cite{Champ} and \cite{Math}. As the latter authors' work does not
seem to be very well known, we provide here a brief discussion of it
here. 

Let us denote by $\partial^{2}\Gamma$ the
space $\{ (a,b) \in \bgamma \rtimes \bgamma \; | \; a
\not= b\}$. Let $\widetilde{G\Gamma}$ denote the
collection of geodesics in $\ebar$. Note that
$\widetilde{G\Gamma}$ has a natural metric with respect to which it is
quasi-isometric to $\ebar$ and hence to $\Gamma$. Furthermore
$\widetilde{G\Gamma}$ carries commuting free and proper actions of $\R$ and
$\Gamma$, and the action of $\Gamma$ is 
co-compact. It is not in general
true that a pair $(a,b)$ of distinct boundary points of $\Gamma$ are connected by a
unique element up to re-parameterization of $\widetilde{G\Gamma}$. In other
words, it is not quite true that $G\Gamma/\R \cong \partial^{2}\Gamma$, which is what
we would like. This may be remedied as follows. 

One defines an equivalence relation $\sim$ on
$\widetilde{G\Gamma}$ such that $G\Gamma = \widetilde{G\Gamma} / \sim$ is
Hausdorff and in fact with
the Hausdorff metric on equivalence classes is a metric space
quasi-isometric to $\widetilde{G\Gamma}$ with the quotient map
$q:\widetilde{G\Gamma} \to G\Gamma$ providing the quasi-isometry. The relation
$\sim$ is $\Gamma$-equivariant, and $\Gamma$ thus acts on $G\Gamma$ and $q$
is a $\Gamma$-invariant map. The relation $\sim$ is not quite
compatible with the action of $\R$ on $\widetilde{G\Gamma}$, but it is
possible to define a new $\R$ action on $G\Gamma$ 
commuting with the $\Gamma$-action and with the following property: if
$(a,b) \in \partial^{2}\Gamma$, the $\R$ orbits of all the geodesics
in $\widetilde{G\Gamma}$ from $a$ to $b$ are collapsed by the quotient map to a single orbit of
the new action of $\R$ on $G\Gamma$. This enables us to identify $G\Gamma/\R$
with $\partial^{2}\Gamma$. 

We remark that this identification may be
seen in another way. If $r$ is a point of $G\Gamma$, the curve $t \mapsto
g_{t}(r)$, where $g_{t}$ denotes the $\R$-action on $G\Gamma$, is a
quasi-geodesic in $G\Gamma$. If under the identification $G\Gamma/\R
\cong\partial^{2}\Gamma$ the $\R$-orbit of $r$ corresponds to $(a,b)
\in \partial^{2}\Gamma$, then it is
also the case that $\lim_{t\to -\infty} g_{t}(r) = a$ and $\lim_{t\to
  +\infty} g_{t}(r) = b$, where the limits are taken in the Gromov
hyperbolic metric space $G\Gamma$. 

We will only need some of the details of this construction in the
proof of Lemma 30. Apart from this lemma, we will only need the
properties of $G\Gamma$ stated in Theorem 27 below.

\begin{rmk} 
\rm

 We choose this moment to note that the only $\Gamma$-invariant homeomorphism $\bgamma \to \bgamma$
is the identity homeomorphism. For, as is well known, the action of
$\Gamma$ on $\bgamma$ is strongly proximal. If $\phi$ is a $\Gamma$-invariant
homeomorphism of $\bgamma$, by amenability of $\Z$, $\phi$ leaves
invariant some probability measure $\mu$. But then for all $\gamma \in
\Gamma$, $\phi_{*}\gamma_{*}(\mu) = \gamma_{*}\mu$. Choose $a \in
\bgamma$. By strong proximality we can choose a
sequence of $\gamma \in \Gamma$ such $\gamma_{*}(\mu) \to \delta_{a}$
where $\delta_{a}$ denotes point mass at $a$, and the convergence is
wk$^{*}$. It follows $\phi$ fixes $a$. Since $a$ was arbitrary, $\phi$
is the identity map.

\end{rmk}

\begin{thm}

There exists a proper metric space $G\Gamma$ on which $\Gamma$
acts, for which:

 1. $G\Gamma$ has the structure of a locally trivial principal $\R$-bundle over
   $\partial^{2}\Gamma$.

 2. $\Gamma$ acts on $G\Gamma$ freely, properly and co-compactly, and its action commutes with the $\R$ action.

 3. There is a continuous involution $G\Gamma \rightarrow G\Gamma$ denoted $r
   \mapsto \hat{r}$, which commutes with the $\Gamma$ action, and satisfies $g_{t}(\hat{r}) =
   \widehat{g_{-t}r}$ for all $t$, where $g_{t}$ denotes the $\R$
   action.

\end{thm}

\begin{note}

\rm

 Elements of the
space $G\Gamma$ should be thought of as geodesics in 
$\ebar$, and so we shall call them $pseudogeodesics$. The $\R$-orbit
of a pseudogeodesic is determined by a pair of distinct
boundary points $(a,b)$. We will call such a pseudogeodesic a
``pseudogeodesic from $a$ to $b$.'' In such a case, we denote by
$r(-\infty)$ the point $a$, and by $r(+\infty)$ the point $b$. As per
the discussion prior to Remark 26, actually the curve $t \mapsto g_{t}(r)$ is a
quasi-geodesic in $G\Gamma$ viewed as a hyperbolic metric space
quasi-isometric to $\Gamma$, and $a = \lim_{t \to -\infty}
g_{t}(r)$, and similarly for $b$, so this notation is actually quite
suitable. 

\end{note}

\begin{rmk}

\rm

If $\Gamma$
acts properly, isometrically and co-compactly on a $ CAT(-\epsilon)$ space $X$
for $\epsilon > 0$ we may take for our purposes the space $G\Gamma$ to be the space of
actual (parameterized) geodesics in $X$, rendering the lemma superfluous. For
convexity in $CAT(-\epsilon)$ spaces implies that any
two distinct boundary points are joined by a unique geodesic.

\end{rmk}

We will also need the following lemma. 

\begin{lemma}

Let $G\Gamma$ be as in Theorem 27. Then there exists a proper $\Gamma$-equivariant map $G\Gamma \rightarrow \ebar$,
    denoted $r \mapsto r(0)$ and satisfying $\lim _{t \rightarrow
    \infty}\;  g_{t}(r) (0) = r(+\infty)$ and $\lim _{t \rightarrow
    -\infty}\;  g_{t}(r) (0) = r(-\infty),$ where the limits are taken
    in the compactified space $\overline{\ebar}$.

\end{lemma}

\begin{proof}

Fixing a point of $G\Gamma$, the orbit map $\Gamma \to G\Gamma$ is a
quasi-isometry which therefore induces a $\Gamma$-invariant homeomorphism
$\bgamma \to \partial G\Gamma$. We may thus identify these two spaces, and
the identification is independent of the point chosen, since any two
such identifications differ by a $\Gamma$-invariant homeomorphism
$\bgamma \to \bgamma$, and the only such is the identity by Remark 26.

On the other hand, by the universal property of $\ebar$
(see \cite{Co3}),
there exists a proper, continuous $\Gamma$-equivariant map $\alpha: G\Gamma \to
\ebar$. Such a map is necessarily a
quasi-isometry, since the action of $\Gamma$ on $G\Gamma$ is
co-compact. Hence $\alpha$ extends to a $\Gamma$-invariant
homeomorphism $\alpha: \bgamma = \partial G\Gamma \to
\bgamma$. Since it is $\Gamma$-invariant, it must be the identity
map, again by Remark 26.

 Now if $r$ is a pseudogeodesic from $a$ to $b$ where $a$ and $b$ are
 points of $\bgamma$ viewed by our identification as points of
 $\partial G\Gamma$, then $t \mapsto g_{t}(r)$ is a quasi-geodesic in
 $G\Gamma$ and we have
$\lim_{t \to -\infty} g_{t}(r) = a$ and $\lim_{t \to +\infty}
g_{t}(r) = b$. Since $\alpha$ is a quasi-isometry, $t \mapsto \alpha
 (g_{t}(r))$ is a quasi-geodesic in $\ebar$, and we have $\lim_{t \to -\infty}
\alpha (g_{t}r) = a$ and $\lim_{t \to +\infty} \alpha
(g_{t} r) = b$ since $\alpha$ is the identity map on the boundary, and we are done. 

Note from this point onward we shall drop the notation $r \mapsto
\alpha
(r)$, replacing it with $r \mapsto r(0)$ as in the statement of the theorem.

\end{proof}

\begin{rmk}

\rm

 Let $M$ be a compact spin$^{c}$ manifold, so that $C(M)$ is a
 Poincar\'e duality algebra in the sense of Definition 12. The fundamental class $\Delta$
 is obtained by pushing forward the class of the Dirac operator on $M$
 by the diagonal map $M \to M \times
 M$ to a class in $K_{*}(M
 \times M) \cong K^{*}(C(M) \otimes C(M))$. Let $U$ be a tubular neighborhood
of the diagonal of in $M \times M$. There is an inclusion of
$C^{*}$-algebras $C_{0}(U) \to C(M) \otimes C(M)$, and the dual
 element $\Dudelta$ is constructed by pushing forward by
this inclusion the Thom class in $K^{*}(U) \cong K_{*}(C_{0}(U))$ to
 an element of $K^{*}(M\times M) \cong K_{*}(C(M) \otimes C(M))$. In our situation, which is vaguely analogous, there is an
inclusion of $C^{*}$-algebras $$C_{0}(\partial^{2}\Gamma) \rtimes
\Gamma \to A \otimes A,$$and the algebra on the left hand side is
strongly Morita equivalent to a cross product by $\R$, and thus has a
Thom class, namely the generator of the flow, which may similarly be
 pushed forward to a class in $K_{1}(A \otimes A)$ and then to a class in $K_{1}(A
 \otimes A^{\rm op})$ using the isomorphism $A \cong A^{\rm op}$. This is how we
shall construct $\Dudelta$. 

\end{rmk}

\begin{note}

\rm

For the following we will denote by $(a,b) \mapsto r_{a,b}$ a
continuous selection of pseudogeodesic from $a$ to $b$. Such a continuous (but not $\Gamma$-equivariant) selection exists by  
Theorem 27.1 and by paracompactness of $\partial^{2}\Gamma$ (see \cite{Gro}).

\end{note}

Define a right $C_{0}(\partial^{2}\Gamma) \rtimes \Gamma$-valued inner
product on the linear space $C_{c}(G\Gamma)$ by the formula: $$<\xi ,
\eta>_{C_{0}(\partial ^{2}\Gamma) \rtimes \Gamma} ((a,b) , \gamma) =
\int_{\R} \bar{\xi} (g_{t} (r_{a,b})) \eta
(g_{t}\gamma^{-1}(r_{a,b}))dt.$$ Define a right $C_{0}(\partial^{2}\Gamma)
\rtimes \Gamma$-module structure on $C_{c}(G\Gamma)$ by $(\xi \cdot f) (r)
= \xi (r) f(r(-\infty) , r(+\infty))$, $f \in
C_{0}(\partial^{2}\Gamma)$, and $(\xi
\cdot \gamma) (r) = \xi (\gamma r)$, for $\gamma \in \Gamma$. Note
this right module structure is
compatible with the inner product.

\begin{defn}

\rm

Let $E$ denote the completion of $C_{c}(G\Gamma)$ to a right Hilbert
$C_{0}(\partial^{2}\Gamma) \rtimes \Gamma$-module with respect to the above
inner product. 

\end{defn}

\begin{defn}

\rm

Define a left action of $C^{*}(\R)$ on $E$ by the unitary
representation $t \mapsto U_{t}$, where
$(U_{t}\xi) (r) = \xi (g_{-t}(r))$.

\end{defn}

\begin{rmk}

\rm

It follows from the definition that the finite rank operators on $E$
as a $C_{0}(\partial^{2}\Gamma) \rtimes \Gamma$-module are linear
combinations of the operators $$K\xi (r) = \sum_{\gamma \in \Gamma}
\zeta (\gamma^{-1}r) \int_{\R} \overline{\eta (g_{t}r)} \xi
(\gamma^{-1}g_{t}r) dt,$$ where $\zeta$ and $\eta$ are elements of
$E$, which fact we will use in the proof (which we have extracted from \cite{Rie}) of
the following lemma.

\end{rmk}

\begin{lemma}

Every element of $C^{*}(\R)$ acts on $E$ as a compact operator. Therefore $E$ defines a class $[E] \in 
  KK(C^{*}( \R) , \; C_{0}(\partial^{2}\Gamma)
  \rtimes \Gamma).$

\end{lemma}

\begin{proof}

(See \cite{Rie}). As $GX/\Gamma$ is compact, we may find a compact fundamental domain
$F$ for the $\Gamma$ action on $GX$. Choose $\epsilon > 0$. Then we
may choose open sets $U_{i}$ of $GX$ such that $F \subset \cup U_{i}$,
and such that for all $i$, $U_{i} \cap g_{t}(U_{i}) = \emptyset$ for
all $|t| \ge \epsilon$. Choose then (see \cite{Rie}) functions
$\zeta_{i , \epsilon} \in C_{c}(GX)$ such that $\zeta_{i,\epsilon}
\in C_{c}(U_{i}),$ and such that  $$(*) \hspace{2cm} \sum_{\gamma \in \Gamma} \zeta_{i, \epsilon}
(\gamma^{-1}r) \int_{\R} \zeta_{i,\epsilon} (\gamma^{-1}g_{t}r) dt =
1 $$ for all $r \in GX$. Define then operators $K_{\epsilon}$ on $E$
by $$K_{\epsilon}\xi (r) = \sum_{i} \sum_{\gamma \in \Gamma} 
\zeta_{i, \epsilon} (\gamma^{-1}r) \int_{\R} \zeta_{i, \epsilon}
(g_{t}r) \xi (g_{t}\gamma^{-1}r) dt.$$ From Remark 35, each
$K_{\epsilon}$ is a compact operator, and from condition $(*)$ above
and the fact that each $\zeta_{i, \epsilon} (r) \zeta_{i, \epsilon}
(g_{t}r) = 0$ if $|t| \ge \epsilon$ and $r \in GX$, it can easily be seen that for $\varphi \in C^{*}(\R)$, $$\varphi \cdot
K_{\epsilon} \to \varphi$$ in operator norm, as $\epsilon \to 0$. Since
each $\varphi \cdot K_{\epsilon}$ is compact, so is $\varphi$.

\end{proof}

\begin{defn}

Let $[D]  =  [\hat{d}_{\R}] \otimes_{C^{*}(\R)} [E] \in KK^{-1}(\C , C_{0}(\partial^{2}\Gamma) \rtimes \Gamma)$, where $[E]$
denotes the class of the cycle $(E, 0)$ for $KK(C^{*}(\R) , C_{0}(\partial^{2}\Gamma) \rtimes \Gamma))$.

\end{defn}

\begin{rmk}

\rm

It will be useful for later to note the following. By the Stabilization
Theorem (\cite{Ka1}) we may embed $E$ as a direct summand of a
trivial Hilbert $C_{0}(\partial^{2}\Gamma) \rtimes \Gamma$-module
$C_{0}(\partial^{2}\Gamma) \rtimes \Gamma \otimes V$, where $V$ is any separable Hilbert
space. Then the left action $C^{*}(\R)
\to B(E)$ of $C^{*}(\R)$ on $E$ may be composed with the embedding,
yielding a homomorphism $\nu: C^{*}(\R) \to K(C_{0}(\partial^{2}\Gamma) \rtimes \Gamma \otimes V) \cong C_{0}(\partial^{2}\Gamma) \rtimes \Gamma \otimes
K(V)$. $[D]$ then becomes $\nu_{*}([\hat{d}_{\R}])$. Note also that since any
two choices of $\nu$ are related by a unitary equivalence, this
construction is not dependent on the choice of embedding $E \to
C_{0}(\partial^{2}\Gamma) \rtimes \Gamma 
\otimes V$.

\end{rmk}

We next note the following trivial:

\begin{lemma} 

The $C^{*}$-algebra $A = C(\bgamma) \rtimes \Gamma$ is isomorphic to
its opposite algebra.

\end{lemma}

\begin{proof}

Define a map $j: A \to A^{\rm op}$ by the covariant pair $j(f) = f$
and $j(\gamma) = \gamma ^{-1}$. Then $j$ induces the required isomorphism.
 
\end{proof}

For what follows, observe that there is a canonical inclusion $ C_{0}(\partial^{2}\Gamma) \rtimes \Gamma \rightarrow A
\otimes A$ given by the composition $C_{0}(\partial^{2}\Gamma) \rtimes
\Gamma \rightarrow C(\bgamma \times \bgamma) \rtimes \Gamma \cong
C(\bgamma) \otimes C(\bgamma) \; \rtimes \Gamma \to C(\bgamma) \rtimes
\Gamma \otimes C(\bgamma) \rtimes \Gamma = A \otimes A$. We
shall denote this inclusion by $i$.

\begin{defn}
\rm

We define the element $\Dudelta \in KK^{-1}(\C , A \otimes A^{\rm op})$ to be $$\Dudelta = (1_{A} \otimes j)_{*}i_{*}([D]) \in KK^{-1}(\C , A \otimes A^{\rm op}).$$
\end{defn}

We are finally in a position to state our main theorem. 

\begin{thm} Let $\Gamma$ be a torsion-free hyperbolic group and $\bgamma$ its
Gromov boundary. Assume that $\bgamma$ has a self-map with no fixed points. Let $A$ denote the cross product $C(\bgamma) \rtimes
\Gamma$. Let $\Delta$ and $\Dudelta$ be the classes constructed in
respectively Definitions 23 and 40. Then $A$ is a Poincar\'e duality algebra in the sense of
Definition 12 and $(\Dudelta , \Delta)$ is a Poincar\'e duality pair.

\end{thm}

The rest of this paper is devoted to the proof of Theorem 41.

\section{Various Reductions}

Let $\Gamma$ be a torsion-free hyperbolic group as in the previous
section, $A$ the cross product $C(\bgamma) \rtimes \Gamma$, and $\Delta \in KK^{1}(A
\otimes A^{\rm op} , \C)$ and $\Dudelta \in KK^{-1}(\C , A \otimes
A^{\rm op})$ the $KK$-classes specified in  respectively Definition 23
and Definition 40. To prove Theorem 41 we must verify that
$$\Dudelta \otimes_{A^{\rm op}}\Delta = 1_{A}$$ and $$\Dudelta
\otimes_{A}\Delta = -1_{A^{\rm op}}.$$ Set $\gamma_{A} = \Dudelta \otimes_{A^{\rm op}} \Delta$ and $\gamma_{A^{\rm
    op}} = \Dudelta \otimes_{A} \Delta$. Using the map $j$ of Lemma 39
we may identify $KK(A^{\rm op} ,
A^{\rm op})$ with $KK(A,A)$. We will first prove that with this
identification, $\gamma_{A}$ and $\gamma_{A^{\rm op}}$ are the same up
to sign, which implies we will only need to compute one of the above products.

 Let $\Delta_{0} = (1_{A} \otimes j)^{*}(\Delta) \in
KK^{1}(A \otimes A , \C)$ and
$\Dudelta_{0} = (1_{A} \otimes j^{-1})_{*}(\Dudelta) \in KK^{-1}(\C ,
A \otimes A)$. Recall that $\sigma_{12}: A \otimes A \to A
\otimes A$ denotes the flip. We first note:

\begin{lemma}

The classes $\Delta_{0}$ and $\Dudelta_{0}$ satisfy $\sigma_{12}^{*}(\Delta_{0}) = \Delta_{0}$; and $(\sigma_{12})_{*}(\Dudelta_{0}) = -\Dudelta_{0}$. 

\end{lemma}

\begin{proof}

Beginning with
$\Dudelta_{0}$, note $\Dudelta_{0} = i_{*}([D])$. Hence it suffices to show $( \sigma _{12} \circ i)_{*} ([D])
= -i _{*}([D])$. Recall $[D]$ is given by $[\hat{d}_{\R}]
\otimes_{C^{*}(\R)} [E]$. Hence $(\sigma_{12} \circ i)_{*}([D]) =
[\hat{d}_{\R}]\otimes_{C^{*}(\R)} (\sigma_{12} \circ i)_{*}[E]$. Let
$u:C^{*}(\R) \to C^{*}(\R)$ denote the homomorphism corresponding to $t \mapsto -t$. Based on a simple index
calculation we see
$u_{*}([\hat{d}_{\R}]) = - [\hat{d}_{\R}]$. Furthermore we have $(\sigma_{12} \circ
i)_{*}([E]) = u^{*}i_{*}([E])$. Hence $(\sigma_{12} \circ i)_{*}([D])
= u_{*}([\hat{d}_{\R}]) \otimes_{C^{*}(\R)} i_{*}([E]) =
-[\hat{d}_{\R}]\otimes_{C^{*}(\R)} i_{*}([E]) = -i_{*}([D])$,
and we are done.

The class $\Delta_{0}$ is represented by the map $A \otimes A \to
Q(H)$, $a\otimes b \mapsto \lambda (a) \rho (b)$, where $\lambda$ is
as before, and $\rho (b) = I\lambda (a) I$, with $I$ the unitary $H
\to H$ induced from inversion on the group. Applying the flip
$\sigma^{*}_{12}$ to $\Delta_{0}$ results in the map $A
\otimes A \to Q(H)$ given by $a\otimes b \mapsto \rho (b)\lambda
(a)$. Since this is conjugate, via $I$, to $\Delta$, the class of these two
extensions is the same: $\sigma_{12}^{*}(\Delta_{0}) = \Delta_{0}$.

\end{proof}

\begin{cor}

We have:  $(j^{-1}_{*})j^{*}(\gamma_{A^{\rm op}}) = -
\gamma_{A}$. Hence if $\gamma_{A} = 1$ then $\gamma_{A^{\rm op}} =
-1_{A^{\rm op}}$. 

\end{cor}

\begin{proof} One checks first that:

\begin{equation} (j^{-1})_{*}(1_{A^{\rm op}}\otimes \Delta) = (j^{-1} \otimes 1_{A}
  \otimes j^{-1})^{*}(1_{A} \otimes \Delta_{0})\end{equation}

 \begin{equation} j^{*}(j^{-1}\otimes 1_{A} \otimes
    j^{-1})_{*}\bigl((\sigma_{12})_{*}(\Dudelta) \otimes 1_{A^{\rm op}}\bigr) =
    (\sigma_{12})_{*}(\Dudelta_{0}) \otimes 1_{A} \end{equation} 

\begin{equation} \gamma_{A^{\rm op}} = \bigl((\sigma_{12})_{*}(\Dudelta)
    \otimes 1_{A^{\rm op}}\bigr) \otimes_{A^{\rm op} \otimes A \otimes
    A^{\rm op}}(1_{A^{\rm op}} \otimes \Delta) .\end{equation}

Hence, using $(3)$, then $(1)$, and then functoriality of the
  intersection product, we have \begin{equation}(j^{-1})_{*}j^{*}( \gamma_{A^{\rm op}}) =
  j^{*}(j^{-1} \otimes 1_{A} \otimes j^{-1})_{*}
 \bigl( (\sigma_{12})_{*}(\Dudelta) \otimes 1_{A^{\rm op}}\bigr)
  \otimes_{A \otimes A \otimes A} (1_{A}
  \otimes \Delta_{0}). \end{equation}

Using $(2)$ we have\begin{equation} (j^{-1})_{*}j^{*}(
  \gamma_{A^{\rm op}}) = \bigl((\sigma_{12})_{*}(\Dudelta_{0}) \otimes
  1_{A}) \otimes (1_{A} \otimes \Delta_{0})\end{equation} 

On the other hand,  \begin{equation} \gamma_{A} = (\Dudelta_{0}
  \otimes 1_{A}) \otimes_{A \otimes A \otimes A} (1_{A} \otimes
  \sigma_{12}^{*}\Delta_{0}), \end{equation} and now, comparing $(5)$ and
  $(6)$ we are done by Lemma 42.

\end{proof}

We are therefore reduced in the proof of Theorem 41 to proving
$\gamma_{A}= 1_{A},$ where, as stated above, $\gamma_{A}$ is the class
$\Dudelta \otimes_{A^{\rm op}}\Delta$. 

\begin{note}

\rm

Recall that if $\mathcal{E}$ is a Hilbert $A$-module, we are denoting
by $B(\mathcal{E})$ the bounded operators on $\mathcal{E}$, $K(\mathcal{E})$ the compact
operators, and $Q(\mathcal{E})$ the quotient
$B(\mathcal{E})/K(\mathcal{E})$. With $\mathcal{E} =
A \otimes H$ the standard Hilbert $A$-module, we have natural maps $A \otimes B(H) \to B(A \otimes H)$, $A
\otimes K(H) \to K(A \otimes H)$ and $A \otimes Q(H) \to Q(A \otimes
H)$. We will sometimes suppress these maps, writing for instance an
element of $B(A \otimes H)$ in the form $a \otimes T$, for $a \in A$
and $T \in B(H)$. 

\end{note}

\begin{rmk}

\rm

For what follows it will be useful to note that any
    function $f$ on $\bgamma \times \Gamma$ continuous in the
    $\bgamma$-variable may be regarded via the formula $f(a)(e_{x}) =
    f(a,x)e_{x}$ as an element of $C(\bgamma ,
    B(H)) \cong C(\bgamma) \otimes B(H)$ whence (see note above), as an element of $B(A
    \otimes H)$, and then, by application of the quotient map, an
    element of $Q(A \otimes H)$.  

\end{rmk}

For further convenience, let us denote the $C^{*}$-algebra
    $C_{0}(\partial^{2}\Gamma) \rtimes \Gamma$ by $B$.

Now, from Equation $(6)$ in the proof of Corollary 43, from
    $\Dudelta_{0} = i_{*}([D])$, and by functoriality of the intersection product, we have $$\gamma_{A}  =
    ([D] \otimes 1_{A})
  \otimes_{ B \otimes A} (i\otimes 1_{A})^{*}\bigl( 1_{A} \otimes
    \sigma_{12}^{*}\Delta_{0}\bigr).$$We will begin by examining the
    term $(i\otimes 1_{A})^{*}\bigl( 1_{A} \otimes
    \sigma_{12}^{*}\Delta_{0}\bigr) \in
  KK^{1}(B\otimes A , A)$. 

Define a covariant pair for the dynamical
  system $(C_{0}(\partial^{2}\Gamma) , \Gamma)$ as follows. If $F$ is
a function on $\partial^{2}\Gamma$  and $\tilde{F}$ denotes
an extension of $F$ to a continuous function on $\bgamma \times
\bar{\Gamma}$, let $\tau(F)$ be the element of
$ Q(A \otimes H)$ corresponding (see Remark 45) to the function
  $\tau(F)(a,x) = \tilde{F}(x^{-1}(a) , x^{-1})$ on $\bgamma \times
  \Gamma$. This is independent of the extension $\tilde{F}$ of $F$. For $\gamma \in \Gamma$, set 
$\tau (\gamma ) =  1 \otimes \lambda^{\rm op}(\gamma ^{-1}) \in Q(A
  \otimes H).$ It is easy to check that these two assignments defines
  a covariant pair.

\begin{defn}

\rm

Let $\tau: B \to Q(A \otimes H)$ be the homomorphism corresponding to the above
covariant pair.

\end{defn}

 For $\gamma \in \Gamma$
recall that $u_{\gamma}$ denotes left translation by $\gamma$. Define
a covariant pair for the dynamical system $(C(\bgamma) , \Gamma)$ by
$\varphi (f) = f \otimes 1 \in B(A \otimes H)$, and $\varphi (\gamma) = \gamma \otimes
 u_{\gamma} \in B( A \otimes H)$.

\begin{defn}

\rm

Let $\varphi: A \to B(A \otimes H)$ denote the homomorphism
corresponding to the above covariant pair. 

\end{defn}

 The following proposition, though depending only on a simple property of
 hyperbolic groups, is central to the proof that $\gamma_{A} =
 1_{A}$. It represents a sort of untwisting of the product $\Dudelta
 \otimes_{A^{\rm op}} \Delta$.

\begin{prop}

The class $(i\otimes 1_{A})^{*}\bigl( 1_{A} \otimes
    \sigma_{12}^{*}\Delta_{0}\bigr) \in
  KK^{1}(B\otimes A , A)$ is represented by the homomorphism $\iota: B
  \otimes A \to Q(A \otimes H)$, $$\iota (b \otimes a ) = \tau
  (b)\pi (\varphi (a)),$$ where $\varphi$, $\tau$ are as in
  Definitions 46 and 47.

\end{prop}

We will require the following: 

\begin{lemma} Let $F \in C_{c}(\partial^{2}\Gamma \times \bgamma)$, and
  $\tilde{F}$ an extension of $F$ to a continuous function on $\bgamma
  \times \bar{\Gamma} \times \bar{\Gamma}$. Then the two functions on
  $\bgamma \times \Gamma$  $$(a,x) \mapsto
  \tilde{F}(x^{-1}(a) , x^{-1} , x) $$and $$(a,x) \mapsto
  \tilde{F}(x^{-1}(a) , x^{-1} , a)$$are the same modulo $C_{0}(\bgamma
  \times \Gamma)$.

\end{lemma}

\begin{proof}

Let $F$ be as in the statement of the lemma. Then for some $\epsilon >
0$, $F$ is supported on the set of $(a,b,c) \in \bgamma \times \bgamma
\times \Gamma$ such that
$d_{\bar{\Gamma}}(a,b) \ge \epsilon$. Therefore $F$ can be
extended to a function $\tilde{F}$ supported on those $(a,b ,
c) \in \bgamma \times \bar{\Gamma} \times \bar{\Gamma}$ for which $d_{\bar{\Gamma}}(a , b) \ge \epsilon$. Let $R$
correspond to $\epsilon$ as in Lemma 14. It suffices to show that for $a \in \bgamma$ fixed and
$x_{n}$ a sequence in $\Gamma$ converging to a boundary point $b \in
\bgamma$, the sequence $$\tilde{F}(x_{n}^{-1}(a ), x_{n}^{-1}, x_{n}) -
\tilde{F}(x_{n}^{-1}(a) , x_{n}^{-1} , a)$$converges to $0$ as $n \to
\infty$. Since if $d_{\bar{\Gamma}}(x_{n}^{-1}a , x_{n}^{-1}) <
\epsilon$, then both $\tilde{F}(x_{n}^{-1}(a ), x_{n}^{-1}, x_{n}) = 0 $ and
$\tilde{F}(x_{n}^{-1}(a) , x_{n}^{-1} , a) = 0$, we may assume after
extracting a subsequence if necessary, that $d_{\bar{\Gamma}}(x_{n}^{-1}(a) , x_{n}^{-1}) \ge
\epsilon$ for all $n$. Then by choice of $R$, $d(x_{0} , [x_{n}^{-1} ,
x_{n}^{-1}a)) = d(x_{n} ,
[e,a)) \le R$ for all $n$, where $[e,a)$ denotes any geodesic ray from $e$ to $a$. Hence $x_{n} \rightarrow a$, and the result follows from continuity of $\tilde{F}$ in the third variable.

\end{proof}

\begin{proof} (Of Proposition 48)

Consider the class $(i\otimes 1_{A})^{*}\bigl( 1_{A} \otimes
    \sigma_{12}^{*}\Delta_{0}\bigr)$. It is represented by the homomorphism $B
  \otimes A \to Q(A \otimes H)$$$a_{1} \otimes a_{2} \otimes
  a_{3} \mapsto a_{1}\otimes  \rho (a_{2}) \lambda (a_{3}),$$where we
  have suppressed the inclusion $i: B \to A \otimes A $ so that in the
  above formula $a_{1} \otimes a_{2}$ is regarded as an element of
  $B$. Here $\rho (a) = \lambda^{\rm op}(j(a))$ as in the proof of
    Lemma 42.  Define a unitary map of Hilbert
modules $U: A \otimes H \to A\otimes H$ by the formula $U (a \otimes
e_{x}) = x \cdot a \otimes e_{x}$. Let $\mathrm{Ad}_{U}$ denote
  the inner automorphism of $Q(A \otimes H)$ given by $\pi(T)
  \mapsto\pi ( UTU^{*})$ and let $\iota '$ denote the homomorphism $B \otimes A \to Q(A
  \otimes H)$ $$\iota ' (a_{1} \otimes a_{2} \otimes a_{3}) =  \mathrm{Ad}_{U}\bigl(a_{1}\otimes  \rho(a_{2}) \lambda (a_{3})\bigr).$$  We claim that
  $\iota ' = \iota.$ It is a simple matter to check that $\iota_{|_{B \otimes C_{r}^{*}(\Gamma)}} =
  \iota ' _{|_{B \otimes C_{r}^{*}(\Gamma)}}$, where $B \otimes
  C_{r}^{*}(\Gamma)$ is viewed as a sub-algebra of $B \otimes A$, and
  that for $b \in B$ and $f \in C(\bgamma)$, we have $\iota (b \otimes f) = \tau (b) \pi (f \otimes 1)$ whereas
  $\iota '(b \otimes f) = \tau (b)\bigl( 1 \otimes \lambda (f)\bigr)$. Thus it remains to prove that $\tau (b) \pi \bigl( 1\otimes  \tilde{f} -
  f\otimes 1\bigr) = 0$ in the Calkin algebra $Q(A \otimes H)$
  whenever $b \in B$, $f \in C(\bgamma)$ and $\tilde{f}$ is an
  extension of $f$ to $\bar{\Gamma}$. Since every $b \in B$ is a
    closed linear combination of elements of the form $\gamma \cdot
    F$, with $\gamma \in \Gamma$ and $F \in C_{c}(\partial ^{2}\Gamma)$, without loss of generality $b =
  F \in C_{c}(\partial^{2}\Gamma)$ and the result follows from Lemma 49.

\end{proof}

\begin{cor}

The class $\gamma_{A}$ lies in the range of the descent map $$\lambda:
RKK_{\Gamma}(\bgamma ; \C , \C) \to KK(A,A),$$i.e. there exists
$\gamma_{\bgamma} \in RKK_{\Gamma}(\bgamma ; \C , \C)$ such that
$\lambda (\gamma_{\bgamma}) = \gamma_{A}$.

\end{cor}

\begin{proof}

Regard (see Remark 38) the class $[D] \in KK^{-1}(\C , B)$ as given by a homomorphism $\nu: C^{*}(\R) \to
B\otimes K(V)$ for some separable Hilbert space $V$. It follows that $[D] \otimes 1_{A}$ is
represented by the homomorphism $\nu\otimes 1_{A}: C^{*}(\R) \otimes A
\to B \otimes A \otimes K(V)$. Hence the class $\gamma_{A}$ is
represented by the homomorphism $C^{*}(\R) \otimes A \to Q(A \otimes H
\otimes V)$ given by the composition $$C^{*}(\R) \otimes A
\stackrel{\nu \otimes 1_{A}}{\longrightarrow} B \otimes A \otimes K(V)
\stackrel{\iota\otimes 1_{K(V)}}{\longrightarrow} Q(A \otimes H
\otimes V).$$ Referring to Lemma 7 with $\Lambda$ the trivial group, let $h$ denote this
composition, and put $h'$ equal to the composition $$C^{*}(\R)
\stackrel{\nu}{\longrightarrow} B\otimes K(V) \stackrel{\tau \otimes
  1_{K(V)}}{\longrightarrow} Q(A \otimes H \otimes V),$$and $h''$
the composition $$A \stackrel{1_{A}}{\longrightarrow} B \otimes A \otimes K(V)
  \stackrel{\iota \otimes 1_{K(V)}} {\longrightarrow} Q(A \otimes H
  \otimes V).$$By Proposition 48, $h''$ lifts to a map $A \to B(A \otimes H \otimes
  V)$ by setting $\tilde{h}''(a) = \varphi (a) \otimes
  1_{V}$. Therefore by Lemma 7, $\gamma_{A}$ is
  represented by the cycle $( A \otimes H \otimes V , F + 1)$, where
  $A \otimes H \otimes V$ has the $(A,A)$-bimodule structure which is
  standard on the right and which on the left is given by the homomorphism $a \mapsto
  \varphi (a) \otimes 1_{V}$, and where $F$ is
  any operator on $A \otimes H \otimes V$ such that $\pi (F) = (\tau
  \otimes 1_{K(V)}) ( \nu(\psi))$.

Now, by construction we may take $F$ to be a limit
of finite linear combinations of operators on the Hilbert
$(A,A)$-bimodule $A \otimes H \otimes V$ of the form \begin{equation}f\otimes
e_{x} \otimes v \mapsto (h_{1}\circ x^{-1})h_{3}\, f \otimes
\tilde{h}_{2}(x^{-1}) e_{x} \otimes T(v)\end{equation}where $T$ is compact, and
$h_{1}\otimes h_{2} \otimes h_{3} \in C_{0}(\partial ^{2}\Gamma \times
\bgamma)$, and where $\tilde{h}_{2}$ denotes a lift of $h_{2}$ to a
continuous function on $\bar{\Gamma}$; and also of the right
translation operators \begin{equation}f\otimes e_{x} \otimes v \mapsto f\otimes
e_{x\gamma}\otimes T(v),\end{equation} where $\gamma \in \Gamma$ and
$T$ is compact. Consider the Hilbert $\bigl(C(\bgamma),
C(\bgamma)\bigr)$-bimodule $C(\bgamma) \otimes H \otimes V$. Let $\Gamma$ act on
$C(\bgamma) \otimes H  \otimes V$ by $\gamma (f
\otimes e_{x} \otimes v) = \gamma (f) \otimes e_{\gamma x} \otimes
v$. Then it is easy to check that with this action, $C(\bgamma)
\otimes H \otimes V$ becomes a $\Gamma - \bigl( C(\bgamma) ,
C(\bgamma) \bigr)$-bimodule. Note that the left and right actions of
$C(\bgamma)$ are in fact the
same. From equations $(7)$ and $(8)$ it is clear that
  $F$ is constructed from operators on $A \otimes H \otimes V$ which
  restrict to operators on $C(\bgamma) \otimes H \otimes V$, hence the
  same is true of $F$. Clearly, as an operator on $C(\bgamma) \otimes
  H \otimes V$, $F$ commutes with
the left action of $C(\bgamma)$ on the module, since this action is the same as
the right action. Finally, $F$ commutes mod
compacts with the action of $\Gamma$, since the operators of which $F$
is built all do. Hence the pair
$(C(\bgamma) \otimes H \otimes V , F + 1)$ actually defines a cycle
for the group $RKK_{\Gamma}(\bgamma ; \C , \C)$. Checking the definition of the
descent map (see \cite{Ka1}) it is easy to see 
that the image of this cycle under descent is precisely the cycle
corresponding to $\gamma_{A}$ described in the first paragraph.

\end{proof}

We will use the above corollary to make use of the following consequence of a theorem of
Tu, which we state in a slightly more general context. Let $\Lambda$
denote a discrete group, which for simplicity we assume acts
co-compactly on its classifying space for proper actions,
$\elambda$ (as is the case for torsion-free hyperbolic $\Gamma$). Let $X$ be a compact metrizable space on which $\Lambda$
acts by homeomorphisms. Recall from Section 1 the map $p_{\elambda } ^{*}: RKK_{\Lambda}(X ;
\C , \C) \rightarrow RKK_{\Lambda}(X \times \elambda  ;\C ,
\C)$. Finally, recall that a $C_{0}(\elambda  \times X )$-algebra $D$
is a $C^{*}$-algebra together with a non-degenerate, asymptotically unital
homomorphism $C_{0}(\elambda  \times X) \to \mathcal{Z}(\mathcal{M}(D))$, where
$\mathcal{Z}$ denotes center. $D$ is called a $\Gamma$-$C_{0}(\elambda
\times X )$-algebra if $\Gamma$ acts by automorphisms on $D$ and the
homomorphism $C_{0}(\elambda  \times X) \to
\mathcal{Z}(\mathcal{M}(D))$ is $\Gamma$-equivariant. Note that such $D$ can be in particular
viewed as a $C(X)$
algebra, by means of the map $C(X) \to C_{b}(\ebar \times
X)) \to \mathcal{Z}(\mathcal{M}(D))$. Let us make the following definition.

\begin{defn}
\rm

Let $D$ be a $\Lambda - C_{0}(\ebar \times X)$-algebra. Define a map
$$\sigma_{\elambda , D}: RKK_{\Lambda}(\elambda \times X; \C , \C) \to
\mathcal{R}KK_{\Lambda}(X ; D , D)$$by replacing a cycle  $(H , F)$ by
 the cycle $(H\otimes_{C_{0}(\elambda \times X)}D , F\otimes 1)$.

\end{defn}

The Hilbert $(D, D)$-bimodule structure on $H\otimes_{C_{0}(\elambda
  \times X)}D$ is well-defined as functions in
$C_{0}(\elambda \times X)$ act as central multipliers of $D$.

Next, we quote Tu's theorem (see \cite{Tu}):

\begin{thm}

 Let the action of $\Lambda$ on $X$ be topologically amenable in the
 sense of \cite{Del}. Then there exist a $\Lambda$-$C_{0}(\elambda
 \times X)$-algebra $D$ and elements $\alpha \in
 \mathcal{R}KK_{\Lambda}(X ; C(X) , D)$, and $\beta \in
 \mathcal{R}KK_{\Lambda}(X; D, C(X)$, satisfying
 $\alpha \otimes_{X , D} \beta = 1_{X} \in \mathcal{R}KK_{\Lambda}(X ;
 C(X) , C(X)) =  RKK_{\Lambda}(X ; \C , \C)$, and $\beta \otimes_{X ,C(X)
 } \alpha = 1_{X , D} \in \mathcal{R}KK_{\Lambda}(X ; D , D)$.

\end{thm}

Using Theorem 52 we can define a map $q: RKK_{\Lambda}(\elambda \times X
; \C , \C) \to RKK_{\Lambda}(X ; \C , \C)$ inverse to  $p_{\elambda}^{*}$ as
follows.

\begin{defn}

\rm

For $a \in RKK_{\Lambda}(\elambda \times X ; \C , \C)$,
define $$q (a) = \alpha \otimes_{X , D}\sigma_{\elambda , D}(a) \otimes_{X , D} \beta \in \mathcal{R}KK_{\Lambda}(X ; C(X) ,
C(X)) = RKK_{\Lambda}(X; \C , \C),$$where $\alpha$ and $\beta$ are as
in Theorem 52 and $\sigma_{\elambda , D}$ is as in Definition 51.

\end{defn}

We show that $q$ and
$p_{\elambda}^{*}$ are inverse to each other. Let $\pi_{1}$ and $\pi_{2}$
denote the projections $\elambda \times \elambda  \to \elambda $, and
$\pi_{1}^{*}$, $ \pi_{2}^{*}$ the corresponding homomorphisms
$C_{0}(\elambda ) \to C_{b}(\elambda  \times \elambda )$. It is a direct consequence of the axioms
for $\elambda $ (see \cite{Co3}) that $\pi_{1}$ and $\pi_{2}$ are $\Lambda$-invariantly
homotopic.

\begin{thm}

The map $p_{\elambda }^{*}$ defines a ring
isomorphism$$RKK_{\Lambda}(X ; \C , \C) \rightarrow
RKK_{\Lambda}(X\times \elambda  ;\C , \C)$$with inverse $q$.

\end{thm}

\begin{proof}

Because the proof is simply an $X$-parameterized version of the
 corresponding statement for $X = \rm pt$ we prove the latter for
 simplicity of exposition. From this assumption we have a $\Lambda -
 C_{0}(\elambda)$-algebra $D$, and $\alpha \in KK_{\Lambda}(\C , D)$,
 $\beta \in KK_{\Lambda}(D , \C)$, satisfying $\alpha \otimes_{D} \beta
 = 1_{\C}$ and $\beta \otimes_{\C} \alpha = 1_{D}$. Let $a \in KK_{\Lambda}(\C ,
 \C)$. Then $q(p_{\elambda}^{*}(a)) = \alpha \otimes_{D} \sigma_{\elambda ,
 D}(p_{\elambda}^{*}(a)) \otimes_{D} \beta =  \alpha \otimes_{D}
 \sigma_{D}(a) \otimes_{D} \beta$, as is easy to check. On the other
 hand, by
 commutativity of the external tensor product and the assumption on
 $\alpha$ and $\beta$, $\alpha \otimes_{D}
 \sigma_{D}(a) \otimes_{D} \beta =  \alpha \otimes_{D} \beta
 \otimes_{\C} a = a$. Hence $q(p_{\elambda}^{*}(a)) = a$.

The other composition is slightly more elaborate. Consider, for $b \in RKK_{\Lambda}(\elambda ; \C ,\C)$ $$p_{\elambda}^{*}(q(b)) = p_{\elambda}^{*}(\alpha) \otimes_{\elambda , D}
p_{\elambda}^{*}(\sigma_{\elambda , D}(b)) \otimes_{\elambda ,
  D}p_{\elambda}^{*}(\beta),$$ and in particular the term
$p_{\elambda}^{*}(\alpha) \otimes_{\elambda , D}p_{\elambda}^{*}(\sigma_{\elambda
  , D}(b))$. We claim that this is equal to $b \otimes_{\elambda} p_{\elambda}^{*}(\alpha)$, whereupon we shall be
done. We may assume that $b$ is given
by a pair $(\mathcal{E} , 0)$, where $\mathcal{E}$ is a
$\Gamma - C_{0}(\elambda)-$module, and that $\alpha$ is given
by a pair $(D , M)$ where $D$ is a $C_{0}(\elambda)$-algebra, and $M$ is a
self-adjoint multiplier of $D$. Then the module for the product $p_{\elambda}^{*}(\alpha) \otimes_{\elambda , D}p_{\elambda}^{*}(\sigma_{\elambda
  , D}(b))$ can be written $\mathcal{E} \otimes_{C_{0}(\elambda)} \bigl(
C_{0}(\elambda) \otimes D \bigr)$, where the tensor product is over the
homomorphism $C_{0}(\elambda) \to C(\elambda \times\elambda) \to
\mathcal{M}(C_{0}(\elambda) \otimes D)$, $f \mapsto \pi_{2}^{*}(f)$. The
operator for the Kasparov product is given by multiplication by $M$ in the
$D$-coordinate; note this is well defined, as $M$, being a multiplier
of $D$, commutes with the actions of functions on $D$. 

On the other hand, consider the product $b \otimes_{\elambda}
p_{\elambda}^{*}(\alpha)$. One
calculates the product of modules to be $\mathcal{E}
\otimes_{C_{0}(\elambda)} \bigl( C_{0}(\elambda) \otimes D \bigr)$,
where this time the tensor product is over the homomorphism $f \mapsto
\pi_{1}^{*}(f)$. The operator is again $M$ acting in the
$D$-coordinate. Now, since $\pi_{1}$ and $\pi_{2}$ are
$\Lambda$-equivariantly homotopic, the two modules are homotopic,
through a homotopy in which the action of $M$ remains the same.

More
precisely, the two cycles corresponding to
the Kasparov products $p_{\elambda}^{*}(\alpha) \otimes_{\elambda , D}p_{\elambda}^{*}(\sigma_{\elambda
  , D}(b))$ and $b \otimes_{\elambda} p_{\elambda}^{*}(\alpha)$ are,
as we have indicated, homotopic, whence $p_{\elambda}^{*}(\alpha) \otimes_{\elambda , D}p_{\elambda}^{*}(\sigma_{\elambda
  , D}(b)) = b \otimes_{\elambda} p_{\elambda}^{*}(\alpha)$. This
proves the claim. (See \cite{Ka1} , pg. 179 for the same sort of argument.)

\end{proof}

\begin{cor}

Let $\gamma_{\bgamma}$ be any element of $RKK_{\Gamma}(\bgamma ;
\C , \C)$ such that $\lambda (\gamma_{\bgamma}) = \gamma_{A}$. Then to show $\gamma_{A} = 1_{A}$, and thus that $(\Dudelta , \Delta)$ is a Poincar\'e duality pair, it suffices to show $p_{\ebar}^{*}(\gamma_{\bgamma }) =
1_{\bgamma \times \ebar}$.

\end{cor}

For $p_{\ebar}^{*}$, being a ring isomorphism, takes a multiplicative unit to a
multiplicative unit. 

Fix $\gamma_{\bgamma}$ to be the class of the cycle for
$RKK_{\Gamma}(\bgamma ; \C , \C)$ described in
the proof of Corollary 50. Then by that corollary $\lambda
(\gamma_{\bgamma}) = \gamma_{A}$. Denote the class $p_{\ebar}^{*}(\gamma_{\bgamma}) \in
RKK_{\Gamma}(\ebar \times \bgamma ; \C , \C)$ by $\gamma_{\ebar \times
  \bgamma}$. By Corollary 55 it remains for us to show that $\gamma_{\ebar \times
  \bgamma} = 1_{\ebar \times \bgamma}$.

\section{Alternative Description of $\gamma_{\ebar
    \times \bgamma}$}

We need first consider more closely the element
$\gamma_{\bgamma }$, as its description in Corollary 50 is unsatisfactory
for our purposes, relying as it does on a inexplicit homomorphism
$\nu: C^{*}(\R) \to B \otimes K(V)$. We would like to describe
a cycle corresponding to $\gamma_{\bgamma }$, whence to $\gamma_{\bgamma
  \times \ebar}$, in such a way as to incorporate the bimodule $E$
associated to the space $G\Gamma$ of pseudogeodesics 
in a more explicit way. Actually, it is quite difficult to do this for
$\gamma_{\bgamma}$ because of dilatability issues, but easy to do it for $\gamma_{\bgamma \times
\ebar}$. So we focus on the latter task. In this section we simply state what this new description
of $\gamma_{\bgamma \times \ebar}$ is, constructing a certain geometric
cycle for
$RKK_{\Gamma}(\bgamma \times \ebar ; \C , \C)$ whose class we will
denote by $\gamma '_{\bgamma \times \ebar}$. We can readily show
that $\gamma_{\ebar \times \bgamma}' = 1_{\ebar \times \bgamma}$. In
the last section we will
verify that in fact $\gamma_{\ebar \times \bgamma} = \gamma ' _{\ebar
  \times \bgamma}.$ Taking these two results together, we will thus have proven $\gamma_{\ebar \times
  \bgamma} = 1_{\ebar \times \bgamma}$.

Recall that we are assuming $\bgamma$ has a fixed point-free
map $S$. By compactness of $\bgamma$ there exists $\delta_{0} > 0$ such that $d_{\bar{\Gamma}}(a, S(a)) \ge \delta_{0}$ for
all $a \in \bgamma$. By abuse of notation, we also denote by $S$ the $equivariant$ map $\bgamma 
\times \Gamma \to \bgamma$ defined by $S(a,z) = z(S(z^{-1}a))$.

\begin{lemma}

There exists an $equivariant$ map $\bgamma \times \Gamma \to
G\Gamma$, $(a,z) \mapsto r_{a,z}$, satisfying $r_{a,z}(-\infty) = a$ and $r_{a,z}(+\infty) =
S(a,z)$. 

\end{lemma}

\begin{proof}

For each $(a,b) \in \partial^{2}\Gamma$, let $r_{a,b}$ be a
  pseudogeodesic from $a$ to $b$, such that the map $(a,b) \mapsto
  r_{a,b}$ is continuous (see Note 32). For $a \in \bgamma$, let
  $r_{a}  =  r_{a ,
  S(a)}$. We have $r_{a}(-\infty ) = a$ and
$r_{a}(+\infty) = S(a)$. To construct an equivariant map as required,
we may set $r_{a,z} = z(r_{z^{-1}a})$.

\end{proof}

Recall that $N$ is the parameter of the Rips complex, which we have
fixed throughout.

\begin{lemma}

There exists a continuous function $Q$ on $\bgamma \times \Gamma
\times \bar{\Gamma}$ satisfying the following properties:

(1) $0 \le Q(a,z,x) \le 1$ for all $(a,z,x) \in \bgamma \times \Gamma \times
\bar{\Gamma} $;

(2) $Q$ is invariant under
the triple diagonal action of $\Gamma$ on $\bgamma \times \Gamma
\times \bar{\Gamma}$;

(3) If $x_{n}$ is a sequence in $\Gamma$, $z \in \Gamma$, and $x_{n}
\to S(a,z)$, then for every $w \in B_{N}(z)$, we have
$Q(a,w,x_{n}) \rightarrow 0$.

(4) If $x_{n}$ is a sequence in $\Gamma$, $z \in \Gamma$, and $x_{n} \rightarrow
a$, then for every $w \in B_{N}(z)$, we have $Q(a,w,x_{n}) \rightarrow 1$.

\end{lemma}

\begin{proof}

Let $Q(a,x)$ be a continuous function on $\bgamma \times
\Gamma$ such that $0 \le Q \le 1$,  $Q(a,x) = 1$ for $d_{\bar{\Gamma}}(a,x)
  < \frac{\delta}{2}$, and $Q(a,x) = 0$ for
  $d_{\bar{\Gamma}}(a,x) \ge \delta$, where $\delta $ is to be
  determined later. Let then $Q(a,z,x)$ be the continuous function on $\bgamma \times \Gamma
\times \bar{\Gamma}$ defined by $Q(a,z,x) = Q(z^{-1}a,
z^{-1}x)$ for $z \in \Gamma$. $Q$ is invariant under the
triple diagonal $\Gamma$ action on $\bgamma \times \bar{\Gamma} \times
\Gamma$. We prove the statement $(3)$; the statement $(4)$ is
similar.

We claim that to prove $Q$ has the required property,we may
assume $z=x_{0}$, where recall $x_{0}$ is the identity of the group
$\Gamma$, regarded as a basepoint in $\ebar$. For, assuming the result
for $z = x_{0}$, let $z$ be arbitrary. Let $w$ be such that $d(z,w) \le
N$, and let $x_{n} \rightarrow S(a , z) = zS(z^{-1}(a))$. Then
$z^{-1}x_{n} \rightarrow S(z^{-1}(a))$. Now $d(z^{-1}w , x_{0}) \le
N$. Hence $Q(z^{-1}(a) , z^{-1}w , z^{-1}x_{n}) \rightarrow 0$ by what we
have assumed proved. But
$Q(z^{-1}(a) , z^{-1}w , z^{-1}x_{n}) = Q(a,w,x_{n})$, by equivariance
of $Q$. This proves the claim.

Let $\delta_{0}$ be as in the paragraph preceding Lemma 56, and let
$R_{0}$ correspond to $\delta_{0}$ as in Lemma 14. Thus, for every $a
\in \bgamma$ we have $d(x_{0} , [a,S(a)]) \le R_{0}$. Choose $ R> 2N
+ 2R_{0}$, choose
$\delta$ according to $R$ as per Lemma 14, and then $Q$ in the first
paragraph as
corresponding to $\delta$. The result of these choices is that $Q(a,x) = 0$ unless $d(x_{0}
, [x,a]) \ge R$.  Let then $x_{n} \rightarrow S(a)$ and let $w \in
B_{N}(x_{0})$. Then if $Q(a,w,x_{n}) = Q(w^{-1}(a) , w^{-1}x_{n})$
does not converge to $0$ we may assume after extracting a subsequence
if necessary that for all large $n$, $d(x_{0} ,
[w^{-1}a , w^{-1}x_{n}]) \ge R$. Hence $d(w ,
[x_{n} , a]) \ge R$. Since $x_{n} \rightarrow S(a)$ it follows that $d(w,
[a,S(a)]) \ge \frac{R}{2}$ and hence $d(x_{0} , [a,S(a)]) \ge
\frac{R}{2} - N > R_{0} $, contradicting choice of $R_{0}$.

\end{proof}

Consider the function $Q(a,z,x)$ constructed in Lemma 57. It will
be convenient to view $Q$ as a function on  $\bgamma \times \Gamma
\times \overline{\ebar}$ satisfying the same properties as the original $Q$; this
is easy to arrange, by reproving Lemma 57 with $\bar{\Gamma}$
replaced by $\overline{\ebar}$. Recall the map $G\Gamma \to \ebar$, $r \mapsto
r(0)$ whose existence was proved in Lemma 30. Define a function $\tilde{Q}$ on $\bgamma \times \Gamma
\times G\Gamma$ by the formula $\tilde{Q}(a,z,r) = Q(a,z,r(0))$. Note $\tilde{Q}$ is
$\Gamma$-invariant.

Define a  $C_{0}(\bgamma \times \ebar)$-valued inner product on the
linear space $C_{c}(\bgamma \times \ebar \times \Gamma \times G\Gamma)$ by
the formula $$<\xi ,
\eta> (a, \mu) = \int_{\Gamma} \int_{\R} \bar{\xi}(a , \mu , z ,
g_{t}r_{a,z})dt d\mu (z).$$

Note the above integral in the $z$-variable is simply a finite
sum, as the support of $\mu \in \ebar$ has diameter at most
$N$. Clearly $C_{c}(\bgamma \times \ebar \times \Gamma \times G\Gamma)$
carries left and right actions of $C_{0}( \bgamma \times \ebar)$, and
these two actions agree, and are compatible with the inner product.

\begin{defn}

\rm

Let $\tilde{\mathcal{E}}$ be the Hilbert $\bigl( C_{0}(\bgamma \times
\ebar), C_{0}(\bgamma \times \ebar)\bigr)$-
bimodule obtained by completing 
$C_{c}(\bgamma \times \ebar \times \Gamma \times G\Gamma)$ with respect to
the above inner product.
\end{defn}

\begin{defn}
\rm

Define an operator $\tilde{P}$ on the Hilbert $C_{0}(\bgamma \times
\ebar)$-module $\tilde{\mathcal{E}}$ by
$$(\tilde{P}\xi)(a,\mu, z , r) = \int_{\Gamma}
\tilde{Q}(a,w,r)\xi(a,\mu,w,r)d\mu (w).$$

\end{defn}

\begin{rmk}

\rm

It is possible to view $\tilde{\mathcal{E}}$ as the sections of a field of Hilbert
spaces $\tilde{H}_{(a,\mu)}$ over $\bgamma \times \ebar$, and the operator $\tilde{P}$ as
corresponding to a field
of operators $\tilde{P}_{(a,\mu)}$, in the following manner. For distinct boundary points
$a$ and $b$, let us denote by $[a,b]$ the fiber over $(a,b)$ in the
map $G\Gamma \rightarrow \partial ^{2} \Gamma$ provided by Theorem 27. Note
that $[a,b]$ has a canonical affine structure, and hence there is in
particular a canonical translation invariant measure on it
corresponding to Lebesgue measure on $\R$. Now, if $\mu$ is a point mass
corresponding to a point $z \in \Gamma$, set
$\tilde{H}_{(a,z)} = L^{2}([a,S(a,z)].$ If $\mu$ is an arbitrary
point of $\ebar$, set
$\tilde{H}_{(a,\mu)}$ to be the completion of the linear space of
functions $\Gamma \to \oplus_{z \in \mathrm{supp}(\mu) \subset \Gamma} \tilde{H}_{(a,z)}$ with
respect to the inner product $<\xi , \eta>_{(a,\mu)} =
\int_{\Gamma} <\xi (z) , \eta (z)>_{\tilde{H}_{(a,z)}} d\mu
(z).$ The operator $\tilde{P}$ corresponds to the following field of operators $\{\tilde{P}_{(a,\mu)}\}$. If $\mu$ is a point mass
corresponding to a point $z \in \Gamma$, $\tilde{P}_{(a,z)}$ is given
by pointwise multiplication by $\tilde{Q}(a,z , \cdot)$ in the variable $r \in [a,S(a,z)]$ on 
$\tilde{H}_{(a,z)}$. If $\mu$ is an arbitrary point of $\ebar$, let $\tilde{P}_{(a,\mu)}$
be defined by 
$\tilde{P}_{(a,\mu)}\xi (z)(r) = \int_{\Gamma} \tilde{Q}(a,w,
r)\xi (w)(r)d\mu (w).$

\end{rmk}

\begin{defn}

\rm

Define a homomorphism $C^{*}(\R) \to B(\tilde{\mathcal{E}})$ by the
unitary representation $t \mapsto U_{t}$, where $(U_{t}\xi)(a,
\mu,z,r) = \xi (a,\mu , z, g_{-t}r)$. 

\end{defn}

As the left $C^{*}(\R)$ action so defined commutes with
the $C_{0}(\bgamma \times
\ebar)$ action, we may view
$\tilde{\mathcal{E}}$ as a $\bigl( C^{*}(\R) \otimes C_{0}(\bgamma \times
\ebar), C_{0}(\bgamma \times
\ebar)\bigr)$-bimodule. Next, note that the triple diagonal action of $\Gamma$ on $\bgamma \times
\ebar \times \Gamma \times G\Gamma$ induces an action of $\Gamma$ on
$C_{0}( \bgamma \times
\ebar \times \Gamma \times G\Gamma)$ as linear maps. It is easy to check
that this action is compatible with the inner product
and right action. It will follow from our remarks below that the
homomorphism $C^{*}(\R) \otimes C_{0}(\bgamma \times \ebar) \to
B(\tilde{\mathcal{E}})$ is $\Gamma$-equivariant. Hence
$\tilde{\mathcal{E}}$ is in fact a $\Gamma - \bigl( C^{*}(\R) \otimes
C_{0}(\bgamma \times \ebar), C_{0}(\bgamma \times
\ebar)\bigr)$-bimodule. 

\begin{rmk}

\rm

Note that from the field perspective, the fact that
$\tilde{\mathcal{E}}$ is a $\Gamma - \bigl( C_{0}(\bgamma \times
\ebar) , C_{0}(\bgamma \times \ebar)\bigr)$-bimodule (which it is in
particular, ignoring the left $C^{*}(\R)$-action) may be re-stated
as:  $\gamma \in
\Gamma$ maps $\tilde{H}_{(a,\mu)}$ isometrically onto
$\tilde{H}_{(\gamma a , \gamma \mu)}$. We note also that we can identify $\tilde{H}_{(a,\mu)}$ in a
  $\Gamma$-equivariant fashion with $L^{2}(\R) \otimes
L^{2}_{\mu}(\Gamma)$. Under this identification the
action of $\gamma \in \Gamma$, $\tilde{H}_{(a,z)} \to
\tilde{H}_{(\gamma a , \gamma z)}$ becomes trivial on the
$L^{2}(\R)$ factor, and the usual action on the $L^{2}_{\mu}(\Gamma)$
factor; and the $C^{*}(\R)$ action becomes trivial on the
$L^{2}_{\mu}(\Gamma)$ factor and the regular representation on the $L^{2}(\R)$
factor. All this follows from using the $\Gamma$-equivariant section $(a,z) \mapsto
r_{a,z}$ to identify $\Gamma$-equivariantly each $[a,S(a,z)]$ with $\R$.

\end{rmk}

\begin{lemma}

The following hold:

(1) The map $C^{*}(\R)\otimes C_{0}(\bgamma \times \ebar) \to B(\tilde{\mathcal{E}})$ is
$\Gamma$ equivariant. 

(2) For $\varphi \in C^{*}(\R)$, $f \in C_{0}(\bgamma \times \ebar)$,
$[\varphi f , \tilde{P}] = f [\varphi , \tilde{P}]$ is
a compact operator.

(3) The operator $\tilde{P}$ is $\Gamma$-equivariant: $[\gamma ,
    \tilde{P}] = 0$ for all $\gamma \in \Gamma$.

(4) $\varphi f (\tilde{P}^{2} -
    \tilde{P})$ and $\varphi f (\tilde{P}^{*} - \tilde{P})$ are compact
     for all $\varphi \in C^{*}(\R)$ and $f\in C_{0}(\bgamma \times \ebar).$

\end{lemma}

\begin{proof}

The first statement is clear. Using the field description, it is easy
to see that to prove
the second statement it suffices to prove that for each $(a,\mu) \in
\bgamma \times \ebar$ the
commutators $[\varphi , \tilde{P}_{(a,\mu)}]$ are compact operators on
$\tilde{H}_{(a,\mu)}$, for $\varphi \in C^{*}(\R)$. Under the identification
$\tilde{H}_{(a,z)} \cong L^{2}(\R)$ pointed out in Remark 62, the operators
$\tilde{P}_{(a,z)}$ become multiplication by functions $\chi_{(a,z)}(t)$ which
satisfy $\lim_{t \to -\infty}\chi_{(a,z)}(t) = 1$ and $\lim_{t \to
    +\infty}\chi_{(a,z)}(t) = 0.$ From this it follows immediately
  that for $\varphi \in C^{*}(\R)$,
the commutator $[\varphi , \tilde{P}_{(a,z)}]$ on
  $\tilde{H}_{(a , z)}$ is compact. Indeed, if $\varphi$ is a compactly
  supported function on $\R$, and $\chi$ is a function with $\lim_{t
    \to -\infty} \chi (t) = 1$ and  $\lim_{t
    \to +\infty} \chi (t) = 0$, it is easy to check that the
  commutator of convolution with $\varphi$ and pointwise
  multiplication by $\chi$ is a compact operator on $L^{2}(\R)$. The
  result for the operators $\tilde{P}_{(a,\mu)}$ follows, since each
  $\tilde{P}_{(a,\mu)}$ is a convex combination of the $P_{(a,z)}$. In an
  exactly analogous way one proves that the operators
  $\varphi \bigl( \tilde{P}_{(a,\mu)}^{2} -
  \tilde{P}_{(a,\mu)}\bigr)$, for $\varphi \in C^{*}(\R)$,
  are compact operators on $\tilde{H}_{(a,\mu)}$, which is part of the fourth
  assertion; self-adjointness follows similarly. Equivariance of
  $\tilde{P}$ is a direct consequence of equivariance of
  the function $\tilde{Q}$.

\end{proof}

We have shown the following:

\begin{cor} The pair $(\tilde{\mathcal{E}} , \tilde{P})$ defines a cycle
  for
  the group $RKK_{\Gamma}^{1}(\bgamma \times \ebar ; C^{*}(\R) ,
  \C)$.  

 \end{cor}

\begin{defn} \rm Let $\gamma_{\bgamma \times \ebar}' = p_{\bgamma
  \times \ebar}^{*}([\hat{d}_{\R}]) \otimes_{\bgamma \times \ebar ,
  C^{*}(\R)} [(\tilde{\mathcal{E}} , \tilde{P})] \in RKK_{\Gamma}(\bgamma
  \times \ebar ; \C , \C)$. 

\end{defn}

\begin{prop} We have: $$\gamma_{\bgamma \times \ebar}'  = 1_{\bgamma \times \ebar}.$$

\end{prop}

\begin{proof}

We first deform the cycle corresponding to $\gamma_{\bgamma \times \ebar}'$ as follows. Identifying the field of
Hilbert spaces $(a,\mu) \mapsto \tilde{H}_{(a,\mu)}$ with the field
$(a,\mu)\mapsto L^{2}_{\mu}(\Gamma ;
L^{2}(\R))$ as in Remark 62, form a homotopy of operators
$\tilde{P}^{t}_{(a,\mu)}$ by the formula $$[\tilde{P}^{t}_{(a,\mu)}
\,\xi ](z) = \int_{\Gamma}[ (1-t) \chi_{(a,w)} + t \chi_{(-\infty ,
  0]} ] \xi (w) d\mu (w),$$where the functions $\chi_{a,w}$ are as in
the proof of Lemma 63. 
It is easy to check this formula defines an operator homotopy in
$RKK_{\Gamma}^{1}( \bgamma \times \ebar ; C^{*}(\R) , \C)$, deforming the
cycle corresponding to  $\gamma_{\bgamma \times \ebar}'$ to the cycle given by the same field of Hilbert spaces, but
with the field of operators given on $\tilde{H}_{(a,\mu)}$ by
$\chi_{(-\infty , 0]} \otimes P_{\mu}$. Now, $P_{\mu}$ is a rank one projection
which is in addition $\Gamma$-invariant. Let $\mu \mapsto \xi_{\mu }$
denote a continuous selection of a unit vector in
$L^{2}_{\mu}(\Gamma)$ for which $P_{\mu}\xi_{\mu} = \xi_{\mu}$ and for
which $\gamma \xi_{\mu} = \xi_{\gamma \mu}$, for any $\gamma \in
\Gamma$. We have$$\tilde{H}_{(a,\mu)} = L^{2}(\R) \otimes [
\xi_{\mu}] \;\;  \oplus \; \; L^{2}(\R) \otimes
L^{2}_{\mu}(\Gamma)^{0},$$where $L^{2}_{\mu}(\Gamma)^{0}$ denotes the
functions in $L^{2}_{\mu}(\Gamma)$ with $\mu$-integral $0$, and
$[\xi_{\mu}]$ denotes the one dimensional linear subspace generated by
$\xi_{\mu}$. With respect to this decomposition, the operator
corresponding to our new deformed cycle is simply $\chi_{(-\infty ,
  0]} \otimes 1 \; \oplus \; 0$, and the $C^{*}(\R)$-action is
diagonal. It follows that the deformed cycle is the direct sum of a
degenerate cycle and the cycle given by the $constant$ field of
Hilbert spaces $L^{2}(\R)$, and operators $\chi_{(-\infty , 0]}$, with
the usual $C^{*}(\R)$-action. The class of the latter is $1_{\bgamma
  \times X}$ by a $\bgamma \times \ebar$-parameterized version of Corollary 3, and the class of the former is $0$ in
$RKK$, and so we are done: $\gamma_{\bgamma \times \ebar}' = 1_{\bgamma \times \ebar}$.

\end{proof}

\section{proof that $\gamma_{\ebar \times \bgamma} = \gamma_{\ebar
    \times \bgamma}'$.}

We now pass to proving $\gamma_{\ebar \times \bgamma} = \gamma_{\ebar
    \times \bgamma}'$. Our strategy for doing this is to define an
    element $b \in RKK_{\Gamma}^{-1}(\bgamma \times \ebar ; \C , B)$
    such that $\gamma_{\bgamma \times \ebar} = p_{\bgamma \times
    \ebar}^{*}([D]) \otimes_{\bgamma \times \ebar} b$. We will then
    separately verify that the axioms for a Kasparov product of $ p_{\bgamma \times
    \ebar}^{*}([D])$ and $b$ are satisfied by the cycle for
    $\gamma_{\bgamma \times \ebar}'$ of the previous section, from
    which we will conclude that $\gamma_{\bgamma \times \ebar}' = \gamma_{\bgamma
    \times \ebar}$.

We first recall the homomorphism $\iota: B \otimes A \to Q(A \otimes
H)$, which in Lemma 48 we showed has the form $\iota (b \otimes a) =
\tau (b) \pi (\varphi (a))$, with $\varphi (f) = f \otimes 1 \in B(A
\otimes H)$ and $\varphi (\gamma) = \gamma \otimes u_{\gamma} \in B( A
\otimes H)$. Let $\Gamma$ act on $C(\bgamma) \otimes H$
diagonally. Then it is clear that $\iota$ restricts to a
$\Gamma$-equivariant homomorphism $B \otimes C(\bgamma) \to
Q(C(\bgamma) \otimes H)$ having the form $b \mapsto \tau (b)$, $f
\mapsto \pi (f \otimes 1)$. We denote this latter $\Gamma$-equivariant
homomorphism $B \otimes C(\bgamma) \to
Q(C(\bgamma) \otimes H)$ by
$\iota_{\bgamma}$.

\begin{rmk}

\rm

A great deal of the complication in this part of the argument arises
from the difficulty in representing the class $\gamma_{\bgamma}$ as a
product of two equivariant classes, even whilst knowing that $\gamma_{\bgamma}$
itself is an equivariant class. Specifically, we do not know whether
or not 
the homomorphism $\iota_{\bgamma}$ is dilatable in the
sense of Definition 5. The idea is that this problem will vanish
when inflating everything over $\ebar$. After doing this, the inflated map, which we
will call $\iota_{\bgamma \times \ebar}$, will in fact
become dilatable, and $\gamma_{\bgamma \times \ebar}$ (though not
$\gamma_{\bgamma})$ will become, as we
would like, a product of two equivariant classes, specifically as $p_{\bgamma \times
    \ebar}^{*}([D]) \otimes_{\bgamma \times \ebar} [\iota_{\bgamma
    \times \ebar}]$. The class $b$ mentioned in the first paragraph of
  this section will be simply the dilation of $\iota_{\bgamma \times
    \ebar}$; i.e. $b = [\iota_{\bgamma \times \ebar}]$.

\end{rmk}

Recall the module $E$ of Definition 33. Choose an embedding of $E$ as
a direct summand of a trivial $B$-module $B \otimes V$ for some Hilbert
space $V$, and denote by $\nu$ the homomorphism $C^{*}(\R) \to B
\otimes K(V)$ obtained by the composition $C^{*}(\R) \to K(E) \to K(B
\otimes V) \cong B
\otimes K(V)$. Let $\nu_{\bgamma }$ denote the homomorphism $C^{*}(\R) \otimes
C(\bgamma) \to B \otimes C(\bgamma) \otimes K(V)$ obtained by
tensoring $\nu$ with the identity on $C(\bgamma)$ and re-arranging
factors. Finally, let $\iota_{\bgamma , V}$ denote the homomorphism $B \otimes
C(\bgamma) \otimes K(V) \to Q(C(\bgamma) \otimes H \otimes V)$
obtained by tensoring $\iota_{\bgamma}$ by the identity homomorphism on $K(V)$ and
re-arranging factors. We have essentially already proved the following
lemma (see the proof of Corollary 50 ), but we restate it for the sake of emphasis. Recall the function $\psi \in C^{*}(\R)$ of Section
1.

\begin{lemma} $\gamma_{\bgamma}$ is represented by any cycle of the form $(C(\bgamma) \otimes H
  \otimes V , F+1)$ where $F \in B(C(\bgamma) \otimes H \otimes V)$ is
  any operator for which $\pi (F) = \iota_{\bgamma , V} (\nu_{\bgamma}
  (\psi \otimes 1))$.

\end{lemma}

Now we tensor all the above data with $\ebar$ as follows. Let firstly
$\iota_{\bgamma \times \ebar}$ denote the homomorphism $B \otimes
C_{0}(\bgamma \times \ebar) \to Q(C_{0}(\bgamma \times \ebar)\otimes
H)$ obtained by tensoring $\iota_{\bgamma}$ with the identity on
$C_{0}(\ebar)$ and re-arranging factors. Let $\iota_{\bgamma \times
  \ebar , V}$ denote the homomorphism $B \otimes C_{0}(\bgamma \times
\ebar) \otimes K(V) \to Q( C_{0}(\bgamma \times \ebar) \otimes H
\otimes V)$ obtained by tensoring $\iota_{\bgamma \times \ebar}$ with
the identity on $K(V)$ and re-arranging factors. Finally, let $\nu_{\bgamma
  \times \ebar}$ denote the homomorphism $C^{*}(\R) \otimes
C_{0}(\bgamma \times \ebar) \to B \otimes C_{0}(\bgamma \times \ebar)
\otimes K(V)$ similarly obtained by tensoring with the identity on
$C_{0}(\bgamma \times \ebar)$ and re-arranging factors. Then just as
above we have:

\begin{lemma}

$\gamma_{\bgamma \times \ebar}$ is represented by any cycle of the
form $\bigl(C_{0}(\bgamma \times \ebar) \otimes H \otimes V , G + 1\bigr)$, where
$G$ is any operator on $C_{0}(\bgamma \times \ebar)\otimes H \otimes V$ satisfying
$\pi (G) = \iota_{\bgamma \times \ebar , V} \bigl(\nu_{\bgamma \times
  \ebar} (\psi \otimes 1)\bigr)$. 

\end{lemma}

Now, suppose we knew that $\iota_{\bgamma \times \ebar}$ was
dilatable. Then $\iota_{\bgamma \times \ebar}$ would define a class 
$b = [\iota_{\bgamma \times \ebar}]$ in $RKK^{1}_{\Gamma}(\bgamma \times
\ebar ; B , \C)$, and the class $\gamma_{\bgamma \times \ebar}$
would then factor in the equivariant category as $\gamma_{\bgamma \times \ebar} = p_{\bgamma
  \times \ebar}^{*}([D]) \otimes_{\bgamma \times \ebar , B} b$. For
emphasis, we state this all explicitly as a proposition, leaving the proof, which is
a standard exercise in Kasparov theory, to the reader. 

\begin{prop}

Let $(\mathcal{E} , P)$ is a cycle for $RKK_{\Gamma}^{1}(\bgamma
    \times \ebar ; B , \C)$ for which 
    there exists an isometry $U: C_{0}(\bgamma \times \ebar) \otimes H \to
    \mathcal{E}$ of Hilbert $C_{0}(\bgamma \times \ebar)$-modules such that for every $f \in C_{0}(\bgamma \times \ebar)$
and $b \in B$: $$\pi (U^{*}P\phi (f\otimes b)PU) =
\iota_{\bgamma \times \ebar}(f \otimes  b) ,$$ where $\phi:
    C_{0}(\bgamma \times \ebar)\otimes B \to
    B(\mathcal{E})$ is the left $C_{0}(\bgamma \times \ebar) \otimes
    B$-structure of $\mathcal{E}$. Then

$$\gamma_{\bgamma \times \ebar} = p_{\bgamma \times \ebar}^{*}([D]) \otimes_{\bgamma
  \times \ebar , B} b,$$ where $b$ denotes the class of $(\mathcal{E} , P)$.

\end{prop}

\begin{rmk}

\rm

After constructing such $b$, it will be possible to describe
$\gamma_{\bgamma \times \ebar}$ without mention of the inexplicit
homomorphism $\nu$. For $p_{\bgamma \times \ebar}^{*}([D])$, in
addition to being represented by the homomorphism $\nu_{\bgamma \times
  \ebar}$,  is alternatively 
represented simply by the pair $(C_{0}(\bgamma \times \ebar) \otimes E  ,
0)$. Hence the product $\gamma_{\bgamma \times \ebar} = p_{\bgamma
  \times \ebar}^{*}([D]) \otimes_{\bgamma \times \ebar , B} b$ will be represented by the
cycle $\bigl(E \otimes_{ B} \mathcal{E} , R\bigr)$, where $R$ is a
$P$-connection. This is how we shall show that $\gamma_{\bgamma \times
  \ebar} = \gamma_{\bgamma \times \ebar}'$. We will find a cycle $(\mathcal{E} , P)$ as in the
  hypothesis of the Proposition 70, such that the resulting cycle $\bigl(E
  \otimes_{ B} \mathcal{E} , R\bigr)$ is homotopic to the cycle $(\tilde{E}
  , \tilde{P})$ described in the previous section. Since the latter
  cycle is homotopic to the cycle for $1_{\bgamma \times \ebar}$, 
  we will conclude $\gamma_{\bgamma \times \ebar} = p_{\bgamma \times \ebar}^{*}([D]) \otimes_{\bgamma
  \times \ebar , B} b = \gamma_{\bgamma \times \ebar}' = 1_{\bgamma
  \times \ebar}$.

\end{rmk}

We now set about construction of the cycle $(\mathcal{E} , P)$ and the
    embedding of $C_{0}(\bgamma \times \ebar) \otimes H$ into $\mathcal{E}$ as above.

Define a $C_{0}(\bgamma \times \ebar )$-valued
   inner product on the linear space $C_{c}(\bgamma \times \ebar \times
    \Gamma ; H)$ by the formula $$<\xi , \eta > (a, \mu) = \int_{\Gamma} <\xi (a , \mu , z ) , \eta (a ,
\mu , z )> d\mu (z).$$Note that the integral is a finite sum, as the
   support of $\mu$ has diameter at most $N$, where $N$ is the
   parameter of the Rips complex.

\begin{defn} \rm Let $\mathcal{E}$ be the right Hilbert $C_{0}(\bgamma
    \times \ebar)$-module obtained by completion of $C_{c}(\bgamma \times \ebar \times
    \Gamma ; H)$ with respect to the above inner product. 

\end{defn}

\begin{defn}

\rm

Define an operator $P$ on $\mathcal{E}$ as follows: let$$P
\xi ( a , \mu , z ) (x) = \int_{\Gamma} Q(a,w , x)\xi (a , \mu , w ,
x)d\mu (w).$$

\end{defn}

Once again the integral is a finite sum.

\begin{defn}

\rm

Define a map $\phi: C_{0}(\bgamma \times \ebar) \otimes B \to B(\mathcal{E})$ by the following covariant
pair.  Let $F \in C_{0}(\partial ^{2}\Gamma )$ and $f \in
C_{0}(\bgamma \times \ebar)$. Define then $\bigl(\phi(f \otimes F)
\xi\bigr) (a, \mu , z)(x)  = f(a, \mu) F\bigl(x^{-1}(a) , x^{-1}S(a,z)\bigr)  \xi (a
,\mu , z)(x).$ For $\gamma \in \Gamma$, define $\phi (\gamma )
\xi  (a , \mu , z) (x) = \xi
(a , \mu , z )( x \gamma ).$

\end{defn}

\begin{rmk}

\rm

As we did in the previous section, we can give a somewhat more
intuitive description of the above data in terms of fields. From this point of view,  $\mathcal{E}$ can
be understood as sections of the continuous, equivariant field of Hilbert
spaces $H_{(a, \mu )} = L^{2}_{\mu}(\Gamma
;H)$. Note that for $\mu$ a point mass at a point $z \in \Gamma
\subset \ebar$, $H_{(a,\mu)}$ is simply $ H$.  The homomorphism $\phi$ can be understood as a
field of homomorphisms $\phi_{(a,\mu)}: B \to B(H_{(a,\mu)})$ as follows: first define, for $(a,z) \in \bgamma \times \Gamma$, a
 homomorphism $\phi_{(a,z)}: B \mapsto B(H)$ by
 $\phi_{(a,z)}(F)(x) = F(x^{-1}(a) , x^{-1}(S(a , z))$ and
 $\phi_{(a,z)}(\gamma) = \lambda^{\rm op} (\gamma^{-1})$. Then define, for $( a, \mu) \in \bgamma \times \ebar$, the homomorphism $\phi_{(a, \mu)}: B \to
 B(H_{(a , \mu)})$ by $\phi_{(a , \mu)}(b)(\xi) (z)(x) =
 \phi_{(a,z)}(b)(\xi (z))(x)$. There is a similar description of the
 operator $P$ as a field of operators $P_{(a,\mu)}$: $(P_{(a,\mu)}\xi )(z)(x) = \int_{\Gamma} Q(a,w,x)\xi (z)(x)d\mu (w)$.

\end{rmk}

Next, note that $\Gamma$ acts on $C_{c}(\bgamma \times \ebar \times \Gamma ;
H)$, and the action is compatible with the $\bigl(C_{0}(\bgamma \times
\ebar)\otimes B , C_{0}(\bgamma \times \ebar)\bigr)$-bimodule structure and the inner product. Hence $\mathcal{E}$
has the structure of a $\Gamma - \bigl( C_{0}(\bgamma \times \ebar)
\otimes B , C_{0}(\bgamma \times \ebar) \bigr)$-bimodule. We have furthermore:

\begin{lemma}

If $f \in C_{0}(\bgamma \times \ebar)$ and $b \in B$, then $[P,\phi
(f \otimes b)]$ is compact.

\end{lemma}

\begin{proof}

Let $F \in
 C_{c}(\partial^{2}\Gamma)$ and $f \in C_{0}(\bgamma \times \ebar)$, and fix $(a,\mu) \in \bgamma \times
 \ebar$ and $z \in \mathrm{supp}(\mu)$. Then we have:  $$\bigl(P\phi
( f\otimes F)\bigr)\xi(a, \mu , z)(x) = f(a,\mu)\int_{\Gamma} Q(a,w,x)F(x^{-1}(a) , x^{-1}S(a,w))\xi
(a , \mu , w)(x) d\mu (w)$$and $$\bigl(\phi ( f \otimes F)P \bigr)\xi(a, \mu , z)(x) =
f(a,\mu)F(x^{-1}(a) , x^{-1}S(a,z))\int_{\Gamma} Q(a,w,x)\xi (a, \mu , w)(x)
d\mu (w).$$Let $x \rightarrow \infty$. Note that for any $w \in \;
$supp$(\mu)$ we have $d(z,w) \le N$. Fix such $w$. Now if the scalar
$F(x^{-1}(a) , x^{-1}S(a,w)) - F(x^{-1}(a) , x^{-1}S(a,z))$ does not
converge to $0$, it follows from the fact that $F \in
C_{c}(\partial^{2}\Gamma)$ and the usual argument, that the distance
from $x$ to the geodesic $[S(a , z) , S(a , w)]$ remains bounded, and
hence that either $x \rightarrow S(a , z)$ or $x \rightarrow S(a ,
w)$. But in either case it follows from Lemma 57 and the fact that
 $d(z,w) \le N$ that both $Q(a,z,x)
\rightarrow 0$ and $Q(a,w,x) \rightarrow 0$. We have shown that the difference
$Q(a,w,x)\bigl(F(x^{-1}a , x^{-1}S(a,w)) - F(x^{-1}a ,
x^{-1}S(a,z))\bigr)$ converges to $0$ as $x \rightarrow \infty$ and
 with $z$ and $w$ fixed. It follows this difference converges to $0$
 uniformly in $z$ and $w$, as the latter range over a finite set. From
 this it follows immediately that the difference of the above two
 expressions represents a compact operator on $\mathcal{E}$.

Finally, to show the commutator $[\phi (f \otimes \gamma) ,
P]$ is compact, observe that $$\bigl(\phi (\gamma) P \phi
(\gamma
^{-1})  - P\bigr) \xi (a, \mu , w , x) = \int_{\Gamma} \bigl(Q(a,w,x\gamma) - Q(a , w, x)\bigr)\xi (a, \mu
, w , x) d\mu (w).$$For every $a$ and every $w$ the function $x
\mapsto Q(a,w, x\gamma) - Q(a,w,x)$ lies in $c_{0}(\Gamma)$, since $Q$
is continuous in the $x$-variable. The
result follows immediately.

\end{proof}

The proof of the following lemma follows the same strategy as that of
the previous one, and we omit it.

\begin{lemma}

$\phi (f \otimes b)(P^{2} - P)$ and $\phi (f \otimes b)(P^{*}-P)$
are both compact operators, for all $b \in B$ and $f \in C_{0}(\bgamma \times \ebar)$.

\end{lemma}

We have shown:

\begin{cor} The pair $(\mathcal{E} , P)$ defines a cycle for
  $RKK^{1}_{\Gamma}(\bgamma \times \ebar ; B , \C)$. 

\end{cor}

\begin{defn} \rm  Let $b \in RKK^{1}_{\Gamma}(\bgamma \times \ebar ; B
  , \C)$ denote the class of the cycle $(\mathcal{E} , P)$ above. 

\end{defn}

We next embed $C_{0}(\bgamma \times \ebar) \otimes H$ into
$\mathcal{E}$ as follows.

\begin{defn}

\rm

Define a  map $U: C_{0}(\bgamma \times \ebar ; H) \to \mathcal{E}$ by
the formula $(U \xi) (a,\mu , w) = \xi (a,\mu)$.

\end{defn}

$U$ is clearly an isometric map of $C_{0}(\bgamma \times \ebar
)$-modules.

\begin{rmk}

\rm

 From the field perspective, $U$ consists of the field of
isometries $U_{(a,\mu)}: H \to H_{(a,\mu)}$ sending $\xi$ to the
  constant function $z \mapsto \xi$. Since each $\mu$ is a probability
  measure, $U$ is indeed isometric. 

\end{rmk}

\begin{prop} The hypothesis of Proposition 70 holds for the pair
  $(\mathcal{E} , P)$, and the isometry $U$ above.

\end{prop}

\begin{proof}

For simplicity of exposition we work with fields. From this point of
view it is easy to see that the homomorphism
$\iota_{\bgamma \times \ebar}$ is given by the field of homomorphisms $\{
(\iota_{\bgamma \times \ebar})_{(a,\mu)} : B \to Q(H)$\} over $\bgamma \times \ebar$, with
$(\iota_{\bgamma \times \ebar})_{(a,\mu)}(F)$ the element of $Q(H)$
corresponding to multiplication by the function $x
\mapsto \tilde{F}(x^{-1}a , x^{-1})$, where $\tilde{F}$ is an
extension of $F$ to a continuous function on $\bgamma \times
\bar{\Gamma}$. Secondly, 
$(\iota_{\bgamma \times \ebar})_{(a,\mu)}(\gamma) = \lambda^{\rm
  op}(\gamma^{-1})$. As mentioned above, the isometric module map $U$ becomes the family of isometries $U_{(a,\mu)}: H \rightarrow H_{(a,\mu)}$,
$U_{(a,\mu)}\xi (w) = \xi$ for all $w \in \rm supp \, (\mu)$. Recall
the homomorphisms $\phi_{(a,\mu)}$ defined in the construction of the
cycle corresponding to the class $b$, and the projections
$P_{(a,\mu)}$. We now wish to show that, for any $b \in B$, the
elements
$$T_{b} = \pi (U^{*}_{(a,\mu)}P_{(a,\mu)}\phi_{(a,\mu)}(b)P_{(a,\mu)}U_{(a,\mu)}) -
(\iota_{\bgamma \times \ebar})_{(a,\mu)}(b)$$are zero in the Calkin algebra of $H$. If $b = \gamma
\in \Gamma$,  it is easy to check that $T_{b}$ is the zero operator, and so we can pass to the
case $b = F \in C_{c}(\partial ^{2}\Gamma)$. In this case, a short
calculation shows that $T_{b}$ corresponds to a diagonal operator,
and, moreover, that to show it is $0$ in the Calkin algebra, it is
enough to show that as $x \rightarrow \infty,$ $$\int_{\Gamma}
Q(a,w,x)\tilde{F}(x^{-1}(a) , x^{-1}S(a,w)) d\mu (w) -
\tilde{F}(x^{-1}(a) , x^{-1}) \rightarrow 0,$$ where $\tilde{F}$ is an
extension of $F$ to a continuous function on $\bgamma \times
\bar{\Gamma}$. Firstly, if $x
\rightarrow a,$ then for large enough $x$, $Q(a,w,x) = 1$ for all $
w\in \rm supp \, (\mu),$ and hence the difference between the above
integral and the integral$$\int \tilde{F}(x^{-1}(a) , x^{-1}S(a,w))
-\tilde{F}(x^{-1}(a) , x^{-1}) d\mu (w)$$ converges to $0$ as $x
\rightarrow a$. Considering the latter integral, for every $w$ in the
integrand we certainly have $d_{\bar{\Gamma}}(x^{-1}S_{w}(a) , x^{-1})
\rightarrow 0$ as $x \rightarrow \infty$, else we would have by the
usual argument that for some $w$, the distance from $x$ to the ray
$[e,S_{w}(a))$ remains bounded, which would imply $x \rightarrow
S_{w}(a)$, thus contradicting $x \rightarrow a$ and $a \not=
S_{w}(a)$. Hence for every $w$ in the integrand
$d_{\bar{\Gamma}}(x^{-1}S_{w}(a) , x^{-1}) \rightarrow 0$ as claimed,
and so the integral converges to $0$ by continuity of $\tilde{F}$ in
the second variable. If $x$ does not converge to $a$, it follows that
$\tilde{F}(x^{-1}(a) , x^{-1}) \rightarrow 0$, and we need only show
the integral also converges to $0$. If it does not, for at least one
$w$, say $w_{1}$, $d_{\bar{\Gamma}}(x^{-1}(a) , x^{-1}S_{w_{1}}(a))$
does not converge to $0$, whence $x \rightarrow a$ or
$S_{w_{1}}(a)$. By assumption $x$ does not converge to $a$ so it
converges to $S_{w_{1}}(a)$. But then by Lemma 57, for $all$ $w$ in
the support of $\mu$, $Q(a,w,x) \rightarrow 0$, since for any such
$w$, $d(w,w_{1}) \le N$, and we are done.

\end{proof}

By Proposition 70 we conclude that $\gamma_{\bgamma \times \ebar} =
p_{\bgamma \times \ebar}^{*}([D]) \otimes_{\bgamma \times \ebar , B}
b$. To show that $\gamma_{\bgamma \times \ebar} = \gamma_{\bgamma
  \times \ebar}'$ it therefore suffices to show that also
$\gamma_{\bgamma \times \ebar}' = p_{\bgamma \times \ebar}^{*}([D]) \otimes_{\bgamma \times \ebar , B}
b$, which we will do by verifying that $\gamma_{\bgamma \times
  \ebar}'$ satisfies the axioms for a Kasparov product of
$p_{\bgamma \times \ebar}^{*}([D])$ and $b$. 

Recall from the discussion in Remark 71 that the product
$p_{\bgamma \times \ebar}^{*}([D]) \otimes_{\bgamma \times \ebar , B}
  b$ is given by the cycle $(E \otimes_{B} \mathcal{E} , R)$ where $R$
  is a $P$-connection. We first observe:

\begin{lemma}

$E \otimes_{B} \mathcal{E} \cong \tilde{\mathcal{E}}$ equivariantly,
and under this isomorphism the  $C^{*}(\R)$ action on $E \otimes_{B}\mathcal{E}$ becomes the action of $C^{*}(\R)$ on
 $\tilde{\mathcal{E}}$ defined in Definition 61.

\end{lemma}

\begin{proof}

To see this, we work from the field point of view, whereapon our statement becomes: for every  $(a,\mu) \in \bgamma \times \ebar$, we have $E\otimes_{B}
H_{(a,\mu)} \cong \tilde{H}_{(a,\mu)}$, and that furthermore, under this
isomorphism, the action of $C^{*}(\R)$ on $E \otimes_{B} H_{(a,
  \mu)}$ corresponds to the action of $C^{*}(\R)$ on
$\tilde{H}_{(a,\mu)}$ described in Remark 62.

 The isomorphism is defined
on the dense subset $C_{c}(G\Gamma) \otimes_{B} L^{2}_{\mu}(\Gamma ; \C
\Gamma)$ of $E\otimes_{B} L^{2}_{\mu}(\Gamma ; H)$ by the composition
of linear maps
$$C_{c}(G\Gamma) \otimes_{B} L^{2}_{\mu}(\Gamma ; \C \Gamma) \cong
\bigl(C_{c}(G\Gamma) \otimes_{B} \C \Gamma \bigr) \otimes
L^{2}_{\mu}(\Gamma ) $$$$\cong C_{c}(G\Gamma)
\otimes_{C_{0}(\partial^{2}\Gamma)}  L^{2}_{\mu}(\Gamma ) \rightarrow
L^{2}_{\mu}\bigl(\Gamma ; \oplus_{z \in \rm supp \; \mu} C_{c}([a, S_{z}(a)])\bigr) \rightarrow
\tilde{H}_{(a,\mu)}.$$The penultimate map is induced by
the restriction map $C_{c}(G\Gamma) \to \oplus_{z \in \rm supp(\mu)} \;
C_{c}( [a,S(a, z)]).$ This composition is isometric with respect to
the various Hilbert module norms. The statement regarding the
$C^{*}(\R)$ actions is obvious.  

\end{proof}

\begin{prop}

We have: $ p_{\bgamma \times \ebar}^{*}([D])\otimes_{\bgamma \times \ebar , B} b = \gamma ' _{\bgamma \times X}.$

\end{prop}

\begin{proof}

 We shall prove this by showing that the operator $\tilde{P}$ is a
 $P$-connection. We work with fields. Taking the product pointwise of
 the modules results in the field of modules $ \tilde{H}_{(a,\mu)}$
 by Lemma 83. We
 show that the operator $\tilde{P}_{(a,\mu)}$ described in Remark 60 is a $P_{(a,\mu)}$-connection. Let
 $\xi \in C_{c}(G\Gamma) \subset E$ and $\theta_{\xi}$ denote the operator
 $H \otimes L^{2}_{\mu}(\Gamma) \rightarrow E\otimes_{B}H \otimes
 L^{2}_{\mu}(\Gamma) $, $\eta \mapsto \xi \otimes_{B} \eta$. By
 \cite{Ka1} we need show that the operator $H \otimes
 L^{2}_{\mu}(\Gamma) \to \tilde{H}_{(a,\mu)}$, $$A_{(a,\mu, \xi)}(\eta) =  \tilde{P}_{(a,\mu)}(\xi
 \otimes_{B} \eta) - \xi \otimes_{B}P_{(a, \mu)}(\eta)$$is a compact
 operator, and show as well that an  adjointed version of this equation
 also represents a compact operator. We
 shall show the first; the second is verified analogously. To calculate
 explicitly the operator $A_{(a,\mu , \xi)}$, assume $\eta \in H
 \otimes L^{2}_{\mu}(\Gamma)$ has the simple form $\eta = e_{x}
 \otimes \alpha$ for $\alpha \in L^{2}_{\mu}(\Gamma)$ and $x \in
 \Gamma$. We have
 $$(A_{(a,\mu , \xi)}\eta) (z) (r) =  \int_{\Gamma} (Q(a,w, r(0)) -
   Q(a,w, x))\alpha (w) \xi (x^{-1}(r))d\mu (w)$$and from this it
   is evident that
     it suffices to show that for $x \rightarrow \infty$ and $w \in
     \rm supp(\mu)$, the
     $L^{2}$-norm of the function $h = h(r) = \bigl(Q(a,w , r(0)) -
       Q(a,w , x)\bigr)\xi (x^{-1}(r))$ of $\tilde{H}_{(a,w)} =
       L^{2}([a,S_{w}(a)]) \cong L^{2}(\R)$
         converges to $0$, since this will express $A_{(a,\mu , \xi)}$
     as a norm limit of finite rank operators.

 Choose $\epsilon
         > 0$.  Then there exists $R>0$ such that if $x$ is large
         enough and
         $r(0) \in B_{R}(x)$, then $|Q(a,w, r(0)) - Q(a,w, x) | <
             \epsilon$, by uniform continuity of $Q(a,w ,
             \cdot)$ and the fact that the Gromov compatification of
             $\Gamma$, and also $\ebar$, is `good' (metric balls in
             the word metric become small in the topology of
             $\bar{\Gamma}$ near the boundary.) Also, as $\xi \in C_{c}(G\Gamma),$ there exists some $R$ for which
        $\xi (r) = 0$ unless $r(0) \in B_{R}(x_{0})$.  It
     follows that for $x$ large enough and $r \in [a, S(a, w)]$, either
     $h(r) = 0$ or $|h(r)| < \epsilon \,  |\xi
     (x^{-1}(r)|$. Consequently, for $x$ sufficiently large,
     $\|h\|_{\tilde{H}_{(a,w)}} < \epsilon \, \|\xi\|_{E}$, and we are done.

\end{proof}

\begin{cor}

We have $$\gamma_{\bgamma \times \ebar}' = \gamma_{\bgamma \times
  \ebar} \in RKK_{\Gamma}(\bgamma \times \ebar ; \C , \C)$$and hence
$\gamma_{\ebar \times \bgamma}  = 1_{\bgamma \times \ebar}$.

\end{cor}

This concludes the proof of Theorem 41.

\end{document}